\documentclass[final]{siamltex}

\usepackage{amsmath}
\usepackage{amsfonts}
\usepackage{graphicx}
\usepackage{epstopdf}
\usepackage{url}
\usepackage{subfigure}
\usepackage{multirow}
\usepackage{array}

\newcommand{\Htwo}{{\mathcal{H}_{2}}}

\title{Recycling BiCG with an Application to Model Reduction\footnotemark[1]}

\author{Kapil Ahuja\footnotemark[2]\ , Eric de Sturler\footnotemark[2]\ ,
Serkan Gugercin\footnotemark[2]\ , and Eun R. Chang\footnotemark[2]}

\begin{document}
\maketitle
\renewcommand{\thefootnote}{\fnsymbol{footnote}}
\footnotetext[1]{This material is based upon work supported by the National Science Foundation under Grant No. {NSF-EAR} 0530643, {NSF-DMS} 1025327, and {NSF-DMS} 0645347.}
\footnotetext[2]{Department of Mathematics, Virginia Tech, Blacksburg, VA 24061.}
\renewcommand{\thefootnote}{\arabic{footnote}}

\begin{abstract}
Science and engineering problems frequently require solving a sequence of
dual linear systems. Besides having to store only few Lanczos vectors,
using the BiConjugate Gradient method (BiCG) to solve dual linear systems
has advantages for specific applications. For example, using BiCG to solve
the dual linear systems arising in interpolatory model reduction provides
a backward error formulation in the model reduction framework.
Using BiCG to evaluate bilinear forms -- for example, in quantum Monte
Carlo (QMC) methods for electronic structure calculations -- leads to a
quadratic error bound. Since our focus is on sequences of dual linear systems,
we introduce {\it recycling BiCG}, a BiCG method that recycles two Krylov
subspaces from one pair of dual linear systems to the next pair.
The derivation of recycling BiCG also builds the foundation for developing
recycling variants of other bi-Lanczos based methods, 
such as CGS, BiCGSTAB, QMR, and TFQMR.

We develop an augmented bi-Lanczos algorithm and a modified two-term recurrence to include recycling in the iteration. The recycle spaces are approximate left and right invariant subspaces corresponding to the eigenvalues closest to the origin. These recycle spaces are found by solving a small generalized eigenvalue problem alongside the dual linear systems being solved in the sequence.

We test our algorithm in two application areas. First, we solve a discretized partial differential equation (PDE) of convection-diffusion type. Such a problem provides well-known test cases that are easy to test and analyze further. Second, we use recycling BiCG in the Iterative Rational Krylov Algorithm (IRKA) for interpolatory model reduction. IRKA requires solving sequences of slowly changing dual linear systems. We analyze the generated recycle spaces and show up to 70\% savings in iterations. For our model reduction test problem, we show that solving the problem without recycling leads to (about) a 50\% increase in runtime.
\end{abstract}



\pagestyle{myheadings}
\thispagestyle{plain}
\markboth{AHUJA, DE STURLER, GUGERCIN, CHANG}{RECYCLING BICG}


\section{Introduction}\label{introSec}
We focus on solving the sequence of dual linear systems,
\begin{equation}\label{initialSys}
    A^{(j)}x^{(j)} = b^{(j)}, \quad A^{(j)*}\tilde{x}^{(j)} = \tilde{b}^{(j)},
\end{equation}
where $A^{(j)} \in \mathbb{C}^{n \times n}$ and 
$b^{(j)}, \tilde{b}^{(j)} \in \mathbb{C}^n$
vary with $j$, the matrices $A^{(j)}$ are large and sparse, 
the solution of the dual system is
relevant, and the change from a pair of systems to the next is small.

In several application areas, there are important advantages to solving dual linear
systems using the BiCG algorithm~\cite{fletcher:bicg}.
BiCG has a short-term recurrence, so very few Lanczos vectors have to be stored.
In addition, using BiCG to solve the dual linear systems arising in
interpolatory model reduction provides a backward stable method 
(with respect to the interpolation
conditions) for computing
a reduced order model~\cite{beattie2010isi} (see Section 5.2).
This makes BiCG attractive even for symmetric positive definite (SPD) systems.
Furthermore, in several applications, such as QMC algorithms~\cite{ahuja:qmc},
we need to evaluate bilinear forms of the
type $u^*A^{-1}w$, where $u, w \in \mathbb{C}^n$ and $A$ is non-Hermitian.
Solving dual linear systems for $u$ and $w$ to compute $u^*A^{-1}w$ provides
a quadratic error bound~\cite{strakos2011bilinear}.

Since BiCG is advantageous for solving dual linear systems and we need to solve a sequence of such systems, we focus on Krylov subspace recycling for BiCG. We refer to our recycling BiCG method as RBiCG. In addition, the BiCG algorithm forms the basis  of other popular bi-Lanczos based algorithms like CGS~\cite{sonneveld:CGS}, BiCGSTAB~\cite{vorst:bicgstab}, QMR~\cite{freund1991qmr}, and TFQMR~\cite{freund1993tfqmr}.  Hence, the derivation of RBiCG is also useful for developing recycling variants of these algorithms~\cite{ahuja:MS}.

The convergence of Krylov subspace methods for solving a linear system, to a great extent,
depends on the spectrum of the matrix, and the deflation of eigenvalues close to the origin
usually improves the convergence rate~\cite{morgan:gmresdr, stathopoulos:deflate}.
If the Krylov subspace is augmented with an eigenvector, then the associated eigenvalue
is effectively deflated. Likewise, for BiCG, it can be shown 
that if the dual Krylov subspace,
$K^i(A^{(j)*},\tilde{r}_0)$, is
augmented with left eigenvectors, the corresponding right eigenvectors are
removed from the primal residual (and vice versa if the primal Krylov subspace is
augmented with right eigenvectors) \cite{Stu99c}.
Therefore, while solving a pair of systems, we select approximate 
left- and right invariant subspaces of $A^{(j)}$ (corresponding to small 
eigenvalues in absolute value), and use these to accelerate the solution 
of the next pair of systems. This process is called 
{\em Krylov subspace recycling}, and leads to faster convergence for 
the next pair of systems.

For solving a single linear system, `recycling' has been used in the
GCROT~\cite{sturler:gcrot} and the GMRES-DR~\cite{morgan:gmresdr} algorithms.
For solving a sequence of linear systems, this idea was first proposed
in~\cite{mike:recycle}, where it is applied to the GCROT and the GCRO-DR
algorithms. Recycling techniques are adapted
to short term recurrences in the RMINRES~\cite{shun:base} algorithm; see \cite{Stu2010a}
for an improved version. GCROT as in~\cite{mike:recycle}, GCRO-DR, and RMINRES all focus on
solving a sequence of single systems rather than a sequence of two dual systems, which is the focus here.
For a fixed matrix with multiple right hand sides deflation-based approaches
are proposed in \cite{abdel:eigBiCG,AbdMor07}.
For a comprehensive discussion of recycling algorithms see~\cite{mike:recycle}.

In addition to testing RBiCG for IRKA~\cite{gugercin:base} for interpolatory model reduction,
we test RBiCG for a model convection-diffusion problem.
PDEs of this type are pervasive in science and engineering,
they lead to nonsymmetric matrices for which
BiCG may be well-suited, and they provide well-known
test cases that are easy to reproduce and to analyze further. Convection-diffusion
problems arise, for example, in the Oseen problem (a fixed-point linearization of
the Navier-Stokes equations), in chemically-reacting flows, heat flow in a medium
with transport, and so on. Moreover, any large discretized PDE leads to a potential
model reduction problem, for example, for uncertainty quantification, for optimizing
an engineering process, or indirectly estimating parameters in the model using
measurements. We analyze the generated recycle spaces for both test problems, and
we show up to 70\% reduction in the iteration count. For our model reduction test
problem, using BiCG instead of RBiCG would take approximately $50\%$ more
time to generate the reduced order model. As recycling is not needed for every
pair of linear systems, this means that the improvement in time for those systems where recycling is actually used is substantially larger (see section 6).

To simplify notation, we drop the superscript $j$ in (\ref{initialSys}). At any particular point in the sequence of systems, we refer to $A x = b$ as the primary system and $A^*\tilde{x} = \tilde{b}$ as the dual system. Throughout the paper, $||\cdot||$ refers to the two-norm, $(\cdot, \cdot)$ refers to the standard inner product, and $e_i$ is the $i$-th canonical basis vector. Unless otherwise stated, we refer to the primary system recycle space and the dual system recycle space collectively as the recycle space.

In the next section, we briefly discuss the BiCG algorithm, and in section~\ref{useSpaceSec} we derive the RBiCG algorithm using a previously computed recycle space. How to compute or update such a recycle space efficiently is discussed in section~\ref{computeSpaceSec}. After explaining the basics of interpolatory model reduction, we discuss how RBiCG is applied in IRKA in section~\ref{secModRed}. We present numerical experiments and results in section \ref{resultsSec} and conclusions in Section \ref{conclSec}.

\section{The BiCG algorithm}\label{backSec}
For the primary system, let $x_0$ be the initial guess with residual $r_0 = b - A x_0$.
Krylov subspace methods, in general, find approximate solutions by projection onto the
Krylov subspace associated with $A$ and $r_0$~\cite{vorst:book}. The $i$-th solution
iterate is given by
\begin{align}\label{initialKrylov}
    x_{i} = x_0 + \varrho_{i},
\end{align}
where $\varrho_{i} \in K^{i}(A,r_{0}) \equiv span\{r_{0},\ Ar_{0},\ A^2 r_{0},\ \cdots,\ A^{i-1} r_{0}\}$ is defined by some projection. The BiCG method defines this projection using the Krylov subspace associated with the dual system, leading to two bi-orthogonal bases and a pair of three-term or coupled two-term recurrences. This method is called the bi-Lanczos method~\cite{lanczos:base,fletcher:bicg}. We initialize the Lanczos vectors as follows:
\begin{align*}
    \begin{array}
        [c]{cc}
        v_1 = r_0 / {|| r_0 ||}, & \tilde{v}_1 = \tilde{r}_0 / {|| \tilde{r}_0 ||}.
    \end{array}
\end{align*}
Defining $V_{i} = [v_{1}\ v_{2}\ \ldots\ v_i]$ and ${\tilde V}_{i} = [{\tilde v}_{1}\ {\tilde v}_{2}\ \ldots\ \tilde{v}_i]$, the ($i+1$)-th Lanczos vectors are given by
\begin{eqnarray*}
        \gamma v_{i+1} = A v_i - V_i \tau \perp \tilde{V}_i, & \qquad &
        \tilde{\gamma} \tilde{v}_{i+1} =  A^* \tilde{v}_i - \tilde{V_i}\tilde{\tau}  \perp {V}_i,
\end{eqnarray*}
where the scalars $\gamma$ and $\tilde{\gamma}$ and the vectors $\tau$ and $\tilde{\tau}$ are to be
determined. This bi-orthogonality condition leads to a pair of $3$-term recurrences (see~\cite{saad:book}), so that computation of the ($i+1$)-th Lanczos vectors requires only the $i$-th and the ($i-1$)-th Lanczos vectors. These 3-term recurrences are called the bi-Lanczos relations, and they are defined as follows:
\begin{align*}
    \begin{array}
        [c]{cc}
            AV_{i} = V_{i+1}\underline{T}_{i} =  V_{i}{T}_{i} + t_{i+1,i}v_{i+1}e_i^T, \\
            A^*\tilde{V}_{i} = \tilde{V}_{i+1}\tilde{\underline{T}}_{i} = \tilde{V}_{i}{\tilde{T}}_{i} + \tilde{t}_{i+1,i}\tilde{v}_{i+1}e_i^T, \\
    \end{array}
\end{align*}
where $T_i$, $\tilde{T}_i$ are $i \times i$ tridiagonal matrices, $t_{i+1,i}$ is the last element of the last row of $\underline{T}_i \in \mathbb{C}^{(i+1) \times i}$, and $\tilde{t}_{i+1,i}$ is the last element of the last row of $\tilde{\underline{T}}_i \in \mathbb{C}^{(i+1) \times i}$.

The next step is to find approximate solutions by projection.
To exploit the efficiency of short-term recurrences in the bi-Lanczos algorithm,
we use the bi-orthogonality condition to define the projection.
This leads to a Petrov-Galerkin approach.
Since the columns of $V_i$ form a basis for $K^{i}(A,r_{0})$, we can define $\varrho_{i}$ in (\ref{initialKrylov}) as $\varrho_{i} = V_iy_{i}$, and the bi-orthogonality (or Petrov-Galerkin) condition then implies
\begin{equation*}
    r_i = b - A(x_0 + \varrho_i) = r_0 - A V_i y_i \perp \tilde{V}_{i}.
\end{equation*}
The vector $y_{i}$ is defined by this orthogonality condition. The solution iterate for the dual system, $\tilde{x}_i$, is similarly defined by $\tilde{x}_{i} = \tilde{x}_0 + \tilde{V}_i \tilde{y}_{i}$  and $\tilde{r}_i \perp V_i$. Further simplifications lead to the standard BiCG algorithm (Algorithm 1)~\cite{fletcher:bicg,vorst:book}.

Next, we briefly discuss the breakdown conditions in BiCG and their
remedies~\cite{greenbaum:book, vorst:book}.
The first breakdown happens when, at any step $i$, $\tilde{r}_i^* r_i = 0$.
This is a breakdown in the underlying bi-Lanczos algorithm and is referred to
as a {\it serious breakdown}. There exist so-called look-ahead
strategies~\cite{freund:lookahead,gutknecht:survey} to avoid this breakdown. In addition,
the two-term recurrence for the solution update requires a pivotless LDU decomposition
of the tridiagonal matrix $T_i$, which may not always exist. This breakdown is
referred to as a breakdown of the {\it second kind}, and it can be avoided by
performing the LDU decomposition with $2 \times 2$ block diagonal
elements~\cite{bank1993csb}. The breakdown conditions in RBiCG are the same,
and similar solutions can be applied. Therefore, and for the sake of brevity,
we do not discuss breakdowns for RBiCG separately, and we'll assume henceforth
in our derivations that breakdowns do not occur. Note that extensive
experiments show that BiCG works well, and that breakdowns rarely happen
in practice \cite{saad:book,gutknecht:survey}.

\begin{figure} {\bf Algorithm 1.} {\it BiCG (adapted from~\cite{vorst:book})} \\
    1. Choose initial guesses $x_0$ and $\tilde{x}_0$. Compute $r_0 = b - Ax_0$ and $\tilde{r}_0 = \tilde{b} - A^* \tilde{x}_0$. \\
    2. {\bf if} $(r_0, \tilde{r}_0) = 0$ {\bf then} initialize $\tilde{x}_0$ to a random vector. \\
    3. Set $p_0 = 0$, $\tilde{p}_0 = 0$, and $\beta_0 = 0$.
    Choose {\tt tol} and {\tt max\_itn}. \\
    4. {\bf for} $i = 1 \ldots $ {\tt max\_itn} {\bf do} \\
    $\diamond$ \;  \;  $p_i = r_{i-1} + \beta_{i-1} p_{i-1}$. \\
    $\diamond$ \;  \;  $\tilde{p}_i = \tilde{r}_{i-1} + \bar{\beta}_{i-1} \tilde{p}_{i-1}$. \\
    $\diamond$ \;  \;  $q_i = A p_i$. \\
    $\diamond$ \;  \;  $\tilde{q}_i = A^* \tilde{p}_i$. \\
    $\diamond$ \;  \;  $\alpha_i = (\tilde{r}_{i-1}, r_{i-1}) / (\tilde{p}_i, q_i)$. \\
    $\diamond$ \;  \;  $x_i = x_{i-1} + \alpha_i p_i$. \\
    $\diamond$ \;  \;  $\tilde{x}_i = \tilde{x}_{i-1} + \bar{\alpha}_i \tilde{p}_i$. \\
    $\diamond$ \;  \;  $r_i = r_{i-1} - \alpha_i q_i$. \\
    $\diamond$ \;  \;  $\tilde{r}_i = \tilde{r}_{i-1} - \bar{\alpha}_i \tilde{q}_i$. \\
    $\diamond$ \;  \;  {\bf if} $||r_i|| \leq $ {\tt tol} and $||\tilde{r}_i|| \leq $ {\tt tol} {\bf then break}. \\
    $\diamond$ \;  \;  $\beta_i = (\tilde{r}_i, r_i) / (\tilde{r}_{i-1}, r_{i-1})$. \\
    5. {\bf end for}.
\end{figure}


\section{Recycling BiCG: Using a Recycle Space}\label{useSpaceSec}
In this section, we modify the BiCG algorithm to use a given recycle space. First, we briefly describe the recycling idea used in the GCRO-DR algorithm. After solving the $j$-th primary system in (\ref{initialSys}), GCRO-DR computes the matrices $U,\ C \in \mathbb{C}^{n \times k}$, such that $\mathrm{range}(U)$ is an approximate invariant subspace of $A^{(j)}$, $A^{(j+1)} U = C$ and $C^* C = I$. It then computes an  orthogonal basis for the Krylov subspace $K^{i}\left(\left(I - CC^*\right)A,\left(I - CC^*\right)r_{0}\right)$.  This produces the Arnoldi relation
\begin{align*}
    \begin{array}
     	[l]{l}
    	AV_i = CC^*AV_i + V_{i+1}\underline{H}_{i} \, \Longleftrightarrow \,
    	(I - CC^*) AV_i = V_{i+1}\underline{H}_{i},
    \end{array}
\end{align*}
where $\underline{H}_{i}$ is an $(i+1) \times i$ upper Hessenberg matrix.
GCRO-DR finds the residual-minimizing solution  over the (direct) sum of the
recycle space, $\mathrm{range}(U)$, and the new search space generated, $\mathrm{range}(V_i)$.

In RBiCG, we use the matrix $U$, derived from an approximate right
invariant subspace of $A^{(j)}$, to define the primary system recycle space,
and compute $C = A^{(j+1)}U$. Similarly, we use the matrix $\tilde{U}$,
derived from an approximate left invariant subspace of $A^{(j)}$, to define
the dual system recycle space, and compute $\tilde{C} = A^{(j+1)*} \tilde{U}$.
Instead of $C$ being an orthogonal matrix, $U$ and $\tilde{U}$ are computed
such that $C$ and $\tilde{C}$ are bi-orthogonal; see Section~\ref{secOrthoC}. The number of vectors selected
for recycling is denoted by $k$, and hence, $U$, $\tilde{U}$, $C$,
and $\tilde{C} \in \mathbb{C}^{n \times k}$.
Next, we derive an augmented bi-Lanczos algorithm that computes
bi-orthogonal bases for the primal
and dual Krylov subspaces. The two-term recurrence for the solution update
in RBiCG is derived in Section \ref{updateSec}.


\subsection{The Augmented Bi-Lanczos Algorithm}\label{augBiSec}
The standard bi-Lanczos algorithm computes columns of $V_i$ and $\tilde{V}_i$ such that, in exact arithmetic, $V_{i} \perp_{b} \tilde{V}_{i}$, where $\perp_{b}$ denotes bi-orthogonality; this implies that $\tilde{V}_{i}^*V_{i}$ is a diagonal matrix. Since we recycle spaces $U$ and $\tilde{U}$, the bi-Lanczos algorithm must be modified to compute the columns of $V_i$ and $\tilde{V}_i$ such that either
\begin{equation}\label{biorthOption}
    \left[U\ V_{i}\right] \perp_{b} \left[\tilde{U}\ \tilde{V}_{i}\right]
\end{equation}
or
\begin{equation}\label{biorth}
    \left[C\ V_{i}\right] \perp_{b} \left[\tilde{C}\ \tilde{V}_{i}\right].
\end{equation}
We choose to implement (\ref{biorth}), because it leads to simpler algebra and hence a more efficient algorithm.
It also has the advantage that the RBiCG algorithm has a form similar to the standard BiCG algorithm.
Next, we derive the recurrences that implement (\ref{biorth}), where $C \perp_b \tilde{C}$
has already been satisfied.
The latter relation is easy to implement when computing the recycle space. Indeed, we can compute
$C$ and $\tilde{C}$ such that $\tilde{C}^*C$ is a real, positive, diagonal matrix; see Section 4.3.
As in the BiCG algorithm, we assume $v_1$ and $\tilde{v}_1$ are available from the initial residuals $r_0$ and $\tilde{r}_0$. We make this statement more precise below.
The ($i+1$)-th Lanczos vector for the primary system is computed by
\begin{align} \label{newVEqn}
    \begin{array}
        [l]{l}
	\gamma v_{i+1} = A v_i - V_i \tau - C \rho \perp \left[\tilde{C}\ \tilde{V}_i\right],
    \end{array}
\end{align}
where $\gamma$, $\tau$, and $\rho$ are to be determined.
Combining (\ref{biorth}) and (\ref{newVEqn}), we get the following equations,
\begin{align}\label{tauRhoEqn}
    \begin{array}
        [l]{l}
        \mathcal{D}_c \rho = \tilde{C}^* A v_i, \\
        \mathcal{D}_i \tau = \tilde{V}_i^* A v_i,
    \end{array}
\end{align}
where $\mathcal{D}_i = \tilde{V}_i^* V_i$ and $\mathcal{D}_c = \tilde{C}^* C$ are both diagonal matrices and
$\mathcal{D}_c$ has real, positive coefficients (see Section \ref{secOrthoC}).
As discussed before, a breakdown in the standard BiCG algorithm because of
singular $\mathcal{D}_i$ can be fixed with look-ahead strategies. Assuming breakdowns do
not occur, we can solve for $\tau$ and $\rho$ in (\ref{tauRhoEqn}) and choose a
normalization $\gamma$; substituting these into (\ref{newVEqn}) gives the ($i+1$)-th Lanczos vector.
Because of the bi-orthogonality condition (\ref{biorth}), the full recurrence for
$v_{i+1}$ reduces to a ($3+k$)-term recurrence, where $k$ is the number of columns of
$C$. This implies that the computation of the ($i+1$)-th Lanczos vector requires
the $i$-th and ($i-1$)-th Lanczos vectors and $C$. Similarly, we get a ($3+k$)-term
recurrence for computing the Lanczos vectors for the dual system. We refer to this
pair of ($3+k$)-term recurrences as the augmented bi-Lanczos relations; they are given by
\begin{eqnarray}\label{augBiLanOld}
	(I-C \hat{C}^*)AV_i = V_{i+1}\underline{T}_i, & \qquad &
	(I-\tilde{C}\check{C}^*)A^* \tilde{V}_i = \tilde{V}_{i+1}\tilde{\underline{T}}_i,
\end{eqnarray}
where
\begin{align}\label{eq:CcheckChat}
	\begin{array}[l]{lll}
	\hat{C}=\left[\frac{\tilde{c}_{1}}{c_{1}^*\tilde{c}_{1}}\ \frac{\tilde{c}_{2}}{c_{2}^*\tilde{c}_{2}}\ \cdots\ \frac{\tilde {c}_{k}}{c_{k}^*\tilde{c}_{k}}\right] & = & \tilde{C}\mathcal{D}_c^{-*} =
    \tilde{C}\mathcal{D}_c^{-1}, \\[5pt]
	\check{C}=\left[\frac{c_{1}}{\tilde{c}_{1}^*c_{1}}\ \frac{c_{2}}{\tilde{c}_{2}^*c_{2}}\ \cdots\ \frac{c_{k}}{\tilde{c}_{k}^*c_{k}}\right] & = & {C}\mathcal{D}_c^{-1}.
	\end{array}
\end{align}
Using (\ref{biorth}), we can rewrite (\ref{augBiLanOld}) as
\begin{align}\label{augBiLanNew}
	\begin{array}[l]{lll}
	A_1 V_i = V_{i+1}\underline{T}_i, & \quad \text{where} & \quad A_1 = (I - C \mathcal{D}_c^{-1} \tilde{C}^*)A(I - C \mathcal{D}_c^{-1} \tilde{C}^*), \\[5pt]
	A_1^* \tilde{V}_i = \tilde{V}_{i+1}\tilde{\underline{T}}_i, & \quad \text{where} & \quad A_1^* = (I - \tilde{C} \mathcal{D}_c^{-*} C^*)A^*(I - \tilde{C} \mathcal{D}_c^{-*} C^*),
	\end{array}
\end{align}
since $\tilde{C}^*V_i = 0$ and $C^*\tilde{V}_i = 0$. This new form of the augmented bi-Lanczos relations simplifies the derivation of the recurrence for the RBiCG solution update, because the operators (\ref{augBiLanNew}) are each other's conjugate transpose. Note that the additional orthogonalizations in (\ref{augBiLanNew}) need not be carried out in
an actual algorithm (see Algorithm~\ref{alg:algo_2}).


\subsection{The Solution Update for the Augmented Bi-Lanczos Recurrence}\label{updateSec}
The $i$-th solution update in the RBiCG algorithm becomes
\begin{equation}\label{solUpdate}
    \begin{array}
    [c]{cc}
        x_{i} = x_{0} + U z_{i} + V_{i} y_{i}, & \tilde{x}_{i} = \tilde{x}_{0} + \tilde{U} \tilde{z}_{i} + \tilde{V}_{i} \tilde{y}_{i},
    \end{array}
\end{equation}
With recycling, the bi-orthogonality condition (\ref{biorth}) defines the Petrov-Galerkin condition,
\begin{equation}\label{orth}
    \begin{array}
    [c]{cc}
        r_{i} = r_0 - AUz_i - AV_iy_i \perp \left[\tilde{C}\ \tilde{V}_{i}\right], & \tilde{r}_{i} = \tilde{r}_0 - A^*\tilde{U}\tilde{z}_i - A^* \tilde{V}_i\tilde{y}_i \perp\left[{C}\ {V}_{i}\right].
    \end{array}
\end{equation}
For the remainder of this section, we focus on the primary system. The derivations for the dual system are analogous. The computation of $z_i$ and $y_i$ can be implemented more efficiently than (\ref{orth}) suggests. Defining $\zeta= ||(I - C{\hat C}^{*})r_{0}||$ and $v_{1}\ =\ {\zeta}^{-1}(I - C{\hat C}^{*})r_0$, we get
\begin{equation}\label{zetaEq}
        r_{0}  =  C{\hat C}^{*} r_{0} + \left(I - C{\hat C}^{*}\right) r_{0} =
        \left[  C \ V_{i+1} \right]  \left[
        \begin{array}
            [c]{c}%
            {\hat C}^{*} r_{0}\\
            \zeta e_{1}%
        \end{array}
        \right].
\end{equation}
Using the augmented bi-Lanczos relation (\ref{augBiLanOld}) we get
\begin{equation}\label{UVEq}
    \begin{array}
        [l]{lll}%
        A\left[  U\ \ V_{i}\right]  \left[
        \begin{array}
            [c]{c}%
            z_{i}\\
            y_{i}%
        \end{array} \right]
        & = & \left[C\ \ \ V_{i+1}\right]  \left[
        \begin{array}
            [c]{cc}%
            I & {\hat{C}}^{\ast}AV_{i}\\
            0 & \underline{T}_{i}%
        \end{array}
        \right]  \left[
        \begin{array}
            [c]{c}%
            z_{i}\\
            y_{i}%
        \end{array}
        \right].
    \end{array}
\end{equation}
Substituting (\ref{zetaEq}) and (\ref{UVEq}) in (\ref{orth}) gives
\begin{equation}\label{updatedRmEq}
    \left[
    \begin{array}
        [c]{c}%
        {\tilde C}^{*}\\
        {\tilde V}_{i}^{*}%
    \end{array}
    \right]  \left[  C \ V_{i+1} \right]  \left(
    \begin{array}
        [c]{c}%
        \left[
        \begin{array}
            [c]{c}%
            {\hat C}^{*} r_{0}\\
            \zeta e_{1}%
        \end{array}
        \right]  - \left[
        \begin{array}
            [c]{cc}%
            I & {\hat C}^{*}AV_{i}\\
            0 & \underline{T}_{i}%
        \end{array}
        \right]  \left[
        \begin{array}
            [c]{c}%
            z_{i}\\
            y_{i}%
        \end{array}
        \right]
    \end{array}
    \right)  = 0.
\end{equation}
Using the bi-orthogonality condition (\ref{biorth}) in the above equation we get\footnote{Note that the length of the vector $e_1$ in (\ref{reducedSystem}) is one less than that of $e_1$ in (\ref{updatedRmEq}), although both denote the first canonical basis vector. Also, $T_i$ in (\ref{reducedSystem}) is $\underline{T}_{i}$ without the last row, and hence is an $i \times i$ tridiagonal matrix.}
\begin{equation}\label{reducedSystem}
    \left[
    \begin{array}
        [c]{c}%
        {\hat{C}}^{\ast}r_{0}\\
        \zeta e_{1}%
    \end{array}
    \right]  -\left[
    \begin{array}
        [c]{cc}%
        I & {\hat{C}}^{\ast}AV_{i}\\
        0 & T_{i}%
    \end{array}
    \right]  \left[
    \begin{array}
        [c]{c}%
        z_{i}\\
        y_{i}%
    \end{array}
    \right] = 0.
\end{equation}
Therefore, $y_{i}$ and $z_{i}$ are given by
\begin{align}\label{yAndz}
    \begin{array}
    [l]{lll}
    T_{i} y_{i} & = & \zeta e_{1}, \\
    z_{i} & = & {\hat C}^{*} r_{0} - {\hat C}^{*}A V_{i} y_{i}.
    \end{array}
\end{align}
Substituting (\ref{yAndz}) in (\ref{solUpdate}) leads to the following solution update:
\begin{equation*}\label{finalUpdate}
    x_{i} = x_{0} + U{\hat C}^{*} r_{0} + (I - U {\hat C}^{*} A) V_{i} y_{i},
\end{equation*}
where $y_i$ is obtained from solving $T_{i}y_{i} = \zeta e_{1}$. All computations here are done with matrix-vector products and  $U {\hat C}^{*} A$ is not computed explicitly.

We introduce a slight change of notation to make future derivations simpler. Let $x_{-1}$
and $\tilde{x}_{-1}$ be the initial guesses and $r_{-1} = b - A x_{-1}$ and $\tilde{r}_{-1} = \tilde{b} - A^* \tilde{x}_{-1}$ the corresponding initial residuals. We define
\begin{align}\label{minusXR}
    \begin{array}{lcr}
	    x_0 = x_{-1} + U \hat{C}^*r_{-1}, & \qquad\qquad\qquad & r_0 = (I - C\hat{C}^*)r_{-1}, \\
	    \tilde{x}_0 = \tilde{x}_{-1} + \tilde{U} \check{C}^*\tilde{r}_{-1}, & \qquad\qquad\qquad & \tilde{r}_0 = (I - \tilde{C}\check{C}^*)\tilde{r}_{-1},
    \end{array}
\end{align}
and follow this convention for $x_0$, $\tilde{x}_0$, $r_0$, and $\tilde{r}_0$ for the rest of the paper. Let
\begin{align*}\label{lduEq}
    \begin{array}
    [l]{lll}
    T_i & = & L_i D_i R_i, \\
    G_i & = & (I - U \hat{C}^*A)V_i R_i^{-1}, \\
    \varphi_i & = & \zeta D_i^{-1}L_i^{-1} e_1.
    \end{array}
\end{align*}
As in the standard BiCG algorithm, an LDU decomposition (without pivoting) of $T_i$ might not always exist. We can avoid this breakdown in the same way as done for BiCG (see Section \ref{backSec}). The two-term recurrence for the solution update of the primary system is now given by
\begin{equation*}\label{finalSolComplexForm}
    x_i = x_{i-1} + \varphi_{i,i}G_i e_i \text {\quad for \quad} i \ge 1,
\end{equation*}
where $\varphi_{i,i}$ is the last entry of the vector $\varphi_i$, and $x_0$ is given by (\ref{minusXR}). An analogous update can be derived for the dual system. Note that we never compute any explicit matrix inverse. The matrices under consideration, $D_i$, $L_i$, and $R_i$, are diagonal, lower triangular, and upper triangular respectively.

This two-term recurrence can be simplified such that $T_i$ is not needed explicitly.
To derive further simplifications, we use the operator $A_1$ (instead of $A$) and
follow steps similar to the ones used in the derivation of BiCG~\cite{gutknecht:survey}. Algorithm 2 provides an outline of RBiCG. Some algorithmic
improvements to make the code faster are not given here; see~\cite{ahuja:MS} for further details.

\begin{figure} \label{alg:algo_2}
{\bf Algorithm 2.} {\it RBiCG} \\
    1. Given $U$ and $\tilde{U}$ compute $\check{C}$ and $\hat{C}$ using (\ref{eq:CcheckChat}).  If $U$ and $\tilde{U}$ are not available, then initialize $U$, $\tilde{U}$, $\check{C}$, and $\hat{C}$ to empty matrices. \\
    2. Choose $x_{-1}$, $\tilde{x}_{-1}$ and compute $x_0$, $\tilde{x}_0$, $r_0$, and $\tilde{r}_{0}$ using (\ref{minusXR}). \\
    3. {\bf if} $(r_0, \tilde{r}_0) = 0$ {\bf then} initialize $\tilde{x}_{-1}$ to a random vector. \\
    4. Set $p_0 = 0$, $\tilde{p}_0 = 0$, and $\beta_0 = 0$.
    Choose {\tt tol} and {\tt max\_itn}. \\
    5. {\bf for} $i = 1 \ldots$ {\tt max\_itn} {\bf do} \\
    \begin{tabular}[b]{ll}
    $\diamond$ \;  \;  $p_i = r_{i-1} + \beta_{i-1} p_{i-1}$; &
                       $\tilde{p}_i = \tilde{r}_{i-1} + \bar{\beta}_{i-1} \tilde{p}_{i-1}$ \\
    $\diamond$ \;  \;  $z_i = Ap_i$; &
                       $\tilde{z}_i = A^*\tilde{p}_i$ \\
    $\diamond$ \;  \;  $\zeta_i = \hat{C}^*z_i$; &
                       $\tilde{\zeta}_i = \check{C}^*\tilde{z}_i$; \\
    $\diamond$ \;  \;  $q_i = z_i - C\zeta_i$; &
                       $\tilde{q}_i = \tilde{z}_i - \tilde{C}\tilde{\zeta}_i$ \\
    $\diamond$ \;  \;  $\alpha_i = (\tilde{r}_{i-1}, r_{i-1}) / (\tilde{p}_i, q_i)$; &
                       $\tilde{\alpha}_i = \bar{\alpha}_i$ \\
    $\diamond$ \;  \;  $\zeta_c = \zeta_c + \alpha_i \zeta_i$; &
                       $\tilde{\zeta}_c = \tilde{\zeta}_c + \tilde{\alpha}_i \tilde{\zeta}_i$ \\
    $\diamond$ \;  \;  $x_i = x_{i-1} + \alpha_i p_i$ &
                       $\tilde{x}_i = \tilde{x}_{i-1} + \tilde{\alpha}_i \tilde{p}_i$ \\
    $\diamond$ \;  \;  $r_i = r_{i-1} - \alpha_i q_i$ &
                       $\tilde{r}_i = \tilde{r}_{i-1} - \tilde{\alpha}_i \tilde{q}_i$ \\
    \multicolumn{2}{l}{$\diamond$ \;  \;  {\bf if} $||r_i|| \leq $ {\tt tol} and
       $||\tilde{r}_i|| \leq $ {\tt tol} {\bf then break}} \\
    $\diamond$ \;  \;  $\beta_i = (\tilde{r}_i, r_i) / (\tilde{r}_{i-1}, r_{i-1})$ & \\
    \end{tabular} \\
    6. {\bf end for} \\
    7. $x_i = x_i - U \zeta_c$; \quad
       $\tilde{x}_i = \tilde{x}_i - \tilde{U}\tilde{\zeta}_c$
\end{figure}


\section{Recycling BiCG: Computing a Recycle Space}\label{computeSpaceSec}
We use the matrices $U$ and $\tilde{U}$ to define the primary and dual system recycle spaces.
The recycle space used in solving a linear system is fixed throughout the RBiCG iteration; however,
the basis of the recycle space for the next pair of linear systems is updated periodically using
the bi-Lanczos vectors. We use harmonic Ritz vectors, with respect to the current Krylov subspace,
to approximate left- and right invariant subspaces cheaply.

We use the following definition \cite{sle:jacobi}. Let $S$ be a subspace of $\mathbb{C}^n$. Then $\lambda \in \mathbb{C}$ is a harmonic Ritz value of $A$ and $u \in S \neq 0$ its corresponding harmonic Ritz vector
with respect to the subspace $\mathcal{W} = AS$ if
\begin{equation} \label{harmonicThm}
  (Au - \lambda u) \perp AS.
\end{equation}
In Section \ref{secHarmonic}, we derive a small generalized eigenvalue problem whose solution gives
the desired approximate invariant subspace. The first pair of systems in our sequence of dual linear systems
requires special attention, since there is no recycle space available at the start. We discuss this case in Section \ref{secFirst}. In Section \ref{secOrthoC}, we describe the construction of the bi-orthogonal $C$ and $\tilde{C}$
in (\ref{biorth}) such that $\mathcal{D}_c = \tilde{C}^*C$ has positive real coefficients.
Although, the generalized eigenvalue problem derived in Section \ref{secHarmonic} is of a small dimension,
it would be expensive to set up in a straightforward manner.
We show how to set up the problem efficiently using recurrences in Section \ref{secSimplify}.

\subsection{Computing an Approximate Invariant Subspace}\label{secHarmonic}
We need a sequence of consecutive Lanczos vectors $v_{i}$ and $\tilde{v}_{i}$ and tridiagonal matrices $T_i$ and $\tilde{T}_i$ to build the recycle space. There is a degree of freedom in choosing the scaling of the Lanczos vectors~\cite{greenbaum:book, gutknecht:survey, saad:book}. The following scaling yields $\tilde{T}_i = T_i^*$
(using (\ref{augBiLanNew}) and (\ref{scalingfactor})):
\begin{equation}\label{scalingfactor}
    ||v_{i}|| = 1, \quad \quad (v_i, \tilde{v}_i) = 1.
\end{equation}
Hence, the Lanczos vectors are computed as follows:
\begin{equation*}
    v_{i} = \frac{r_{i-1}}{||r_{i-1}||}, \quad \quad \tilde{v}_{i} = \frac{\tilde{r}_{i-1}}{(v_{i}, \tilde{r}_{i-1})}.
\end{equation*}
$T_i$ can be computed using the residuals and iteration scalars of the
RBiCG iteration as follows~\cite{abdel:eigBiCG,ahuja:MS}:
\begin{align*}
        {T}_{i} =
        \begin{pmatrix}
            \frac{1}{\alpha_1} & -\frac{||r_{0}||}{||r_{1}||} \cdot \frac{\beta_{1}}{\alpha_{1}} &  &  &  \\
            -\frac{||r_{1}||}{||r_{0}||} \cdot \frac{1}{\alpha_1} &  \frac{1}{\alpha_2} + \frac{\beta_{1}}{\alpha_{1}} & -\frac{||r_{1}||}{||r_{2}||} \cdot \frac{\beta_{2}}{\alpha_{2}} &  & \\
            &  & \cdots &  &  \\
            &  & \ddots &  & -\frac{||r_{i-2}||}{||r_{i-1}||} \cdot \frac{\beta_{i-1}}{\alpha_{i-1}} \\
            &  &  &  -\frac{||r_{i-1}||}{||r_{i-2}||} \cdot \frac{1}{\alpha_{i-1}} & \frac{1}{\alpha_i} + \frac{\beta_{i-1}}{\alpha_{i-1}}
        \end{pmatrix}.
\end{align*}

Instead of using all the Lanczos vectors to update the recycle space, we update the recycle space periodically. This strategy keeps the memory requirements modest~\cite{shun:base}, as it allows us to discard Lanczos vectors periodically. The iteration process between two updates of the recycle space is referred to as a ``cycle''. The length of the cycle, $s$, refers to the number of iterations between updates. Let $V_j$ and $\tilde{V}_j$ contain the Lanczos vectors generated during the $j^{th}$ cycle,
\begin{align*}
    \begin{array}
    [l]{ll}
    V_j = \left[v_{(j-1)s + 1} \quad \ldots \quad v_{js}\right], &
    \tilde{V}_j = \left[\tilde{v}_{(j-1)s + 1} \quad \ldots \quad \tilde{v}_{js}\right].
    \end{array}
\end{align*}
Also, let
\begin{align*}
    \begin{array}
    [l]{ll}
    \Upsilon_j = \left[v_{(j-1)s} \quad  V_j \quad v_{js+1}\right], &
    \tilde{\Upsilon}_j = \left[\tilde{v}_{(j-1)s} \quad \tilde{V}_j \quad \tilde{v}_{js+1}\right],
    \end{array}
\end{align*}
where $v_{(j-1)s}$ and $\tilde{v}_{(j-1)s}$ are the last Lanzos vectors from the previous
cycle,  and $v_{js+1}$ and $\tilde{v}_{js+1}$ are the first Lanzos vectors from the
next cycle. The augmented bi-Lanczos relations for the $j^{th}$ cycle are now given by
\begin{eqnarray}\label{augBiLanExtend}
    (I - C\hat{C}^*)AV_{j} = \Upsilon_{j}\Gamma_{j} , & \qquad &
    (I - \tilde{C}\check{C}^*)A^*\tilde{V}_{j} = \tilde{\Upsilon}_{j}\text{ }\tilde{\Gamma}_{j},
\end{eqnarray}
where $\Gamma_{j}$, $\tilde{\Gamma}_{j}  \in \mathbb{C}^{(s+2) \times s}$ are $T_j$, $\tilde{T}_j$, respectively, with an extra row at the top (corresponding to $v_{(j-1)s}$ and $\tilde{v}_{(j-1)s}$) and at the bottom (corresponding to $v_{js+1}$ and $\tilde{v}_{js+1}$).

The discussion in this paragraph concerns only the primary system. However, an analogous discussion applies to the dual system. Let $U$ define the recycle space available from the previous linear system and $U_{j-1}$ the recycle space generated at the end of cycle $(j-1)$ for the current linear system. We want to obtain an improved $U_j$ from $V_j$, $U_{j-1}$, and $U$. It is important to note that $U_j$ is not used for solving the current linear system. At the end of solving the current linear system, the final $U_j$ will be $U$ for the next linear system. There are several choices for selecting $U_j$ \cite{shun:base}. For simplicity, we build $U_j$ from $\mathrm{range}([U_{j-1}\ V_j])$.

Based on the choices discussed in the previous two paragraphs, we first define
certain matrices, and then we derive the generalized eigenvalue problem whose
solution gives the approximate invariant subspace. Let
\begin{align*}
	\begin{array}
	[l]{lll}
	\Phi_{j} =\left[U_{j-1}\ V_{j}\right], &
	\Psi_{j} = \left[C\ C_{j-1}\ \Upsilon_{j}\right], &
	H_{j} = \left[
	\begin{array}
		[c]{cc}%
		0 & B_{j}\\
        		I & 0\\
        		0 & \Gamma_{j}%
   	 \end{array}
    	\right], \\[15pt]
    	\tilde{\Phi}_{j} = \left[\tilde{U}_{j-1}\ \tilde{V}_{j}\right], &
    	\tilde{\Psi}_{j} = \left[\tilde{C}\ \tilde{C}_{j-1}\ \tilde{\Upsilon}_{j}\right], &
    	\tilde{H}_{j} = \left[
    	\begin{array}
        		[c]{cc}%
        		0 & {\tilde B}_{j}\\
        		I & 0\\
        		0 & \tilde{\Gamma}_{j}
    	\end{array}
    	\right],
	\end{array}
\end{align*}
where $C_{j-1} = AU_{j-1}$, $B_{j} = \hat{C}^*AV_{j}$, $\tilde{C}_{j-1} = A^*\tilde{U}_{j-1}$, and ${\tilde B}_{j} = \check{C}^*A^*\tilde {V}_{j}$. Then, the
augmented bi-Lanczos relations (\ref{augBiLanExtend}) lead to
\begin{align}\label{firstKW}
	\begin{array}
	[l]{l}
	A\Phi_{j}=\Psi_{j}H_{j}, \qquad
	A^*\tilde{\Phi}_{j} = \tilde{\Psi}_{j}\tilde {H}_{j}.
	\end{array}
\end{align}
In RMINRES~\cite{shun:base}, harmonic Ritz pairs of $A$ with respect to the subspace $\mathrm{range}(A\Phi_j)$ have been successfully used to build the recycle space. Since we work in a Petrov-Galerkin framework, it is more
intuitive to use harmonic Ritz pairs with respect to the subspace $\mathrm{range}(A^* \tilde{\Phi}_j)$, following~\cite{beattie1998hr}. This leads to simpler algebra and cheaper computations. Let $(\lambda, u)$ denote an harmonic Ritz pair of $A$. Then, we derive $\lambda$ and $u \in \mathrm{range}(\Phi_j)$ from the condition
\begin{equation}\label{harmonic1}
    \left(Au - \lambda u\right) \perp \mathrm{range}\left(A^*\tilde{\Phi}_{j}\right).
\end{equation}
Taking $u = \Phi_{j}w$ and substituting (\ref{firstKW}) in (\ref{harmonic1}) gives
\begin{align*}
    \begin{array}[c]{ccc}
        \left(A^*\tilde{\Phi}_{j}\right)^*A\Phi_{j}w = \lambda\left(A^*\tilde{\Phi}_{j}\right)^*\Phi_{j}w
        & \Leftrightarrow &
        \left(\tilde{\Psi}_{j}\tilde{H}_{j}\right)^*\Psi_{j}H_{j}w = \lambda\left(\tilde{\Psi}_{j}\tilde{H}_{j}\right)^*\Phi_{j}w.
    \end{array}
\end{align*}
Thus, condition (\ref{harmonic1}) leads to the generalized eigenvalue problem,
\begin{equation}\label{finalGenEigValue}
    \tilde{H}_{j}^*\tilde{\Psi}_{j}^*\Psi_{j}H_{j}w = \lambda\tilde{H}_{j}^*\tilde{\Psi}_{j}^*\Phi_{j}w.
\end{equation}
Let the columns of $W_j$ be the $k$ right eigenvectors
corresponding to the eigenvalues closest to the origin.
Then, we take $U_j = \Phi_j W_j$.
See~\cite{ahuja:MS} for an analogous derivation of the dual system recycle space.

\subsection{The First Linear System and the First Cycle}\label{secFirst}
For the first cycle of the first system,
the matrices $U$, $\tilde{U}$, $U_{j-1}$, and $\tilde{U}_{j-1}$ are not available.
Letting $T_1$ and $\tilde{T}_1$ denote the tridiagonal matrices for the first cycle, we
consider the following eigenvalue problems:
\begin{align*}
    T_1 w = \lambda w, \quad \quad \tilde{T}_1 \tilde{w} = \mu \tilde{w}.
\end{align*}
Since $\tilde{T}_1 = T_1^*$, we solve for the left and the right eigenvectors of $T_1$,
$\tilde{W}_1$ and $W_1$ respectively.
Hence, we take
\begin{align*}
    U_1 = V_1 W_1, \quad \quad \tilde{U}_1 = \tilde{V}_1 \tilde{W}_1.
\end{align*}

During the second and subsequent cycles of the first linear system, $U_{j-1}$
and $\tilde{U}_{j-1}$ are available, but $C$ and $\tilde{C}$ are not.
Redefining $\Psi_{j}, \tilde{\Psi}_{j}, H_{j},$ and $\tilde{H}_{j}$, we get
the generalized eigenvalue problem (\ref{finalGenEigValue}) with
\begin{align*}
	\begin{array}
	[l]{lll}
    	\Phi_{j} = \left[U_{j-1} \quad V_j\right], &
	\Psi_{j}=\left[C_{j-1} \quad \Upsilon_{j}\right], &
	H_{j}=\left[
        		\begin{array}
	        [c]{cc}%
	        I & 0\\
	        0 & \Gamma_{j}%
	        \end{array}
	 	\right], \\[10pt]
	 \tilde{\Phi}_{j} = \left[\tilde{U}_{j-1} \quad \tilde{V}_j\right], &
	 \tilde{\Psi}_{j}=\left[\tilde{C}_{j-1} \quad \tilde{\Upsilon}_j\right], &
	 \tilde{H}_{j} = \left[
        		\begin{array}
	        [c]{cc}%
	        I & 0\\
	        0 & \tilde{\Gamma}_{j}
	        \end{array}
	  	\right].
	\end{array}
\end{align*}

For the first cycle of each of the subsequent linear systems (i.e. $j$ = 1), $C$ and $\tilde{C}$ are available, while $U_{j-1}$ and $\tilde{U}_{j-1}$ are not. Redefining $ \Phi_{1} ,  \tilde{\Phi}_{1}, \Psi_{1}, \tilde{\Psi}_{1}, H_{1},$ and $\tilde{H}_{1}$, we get the generalized eigenvalue problem (\ref{finalGenEigValue}) with
\begin{align*}
	\begin{array}
	[l]{lll}
	\Phi_{1} = \left[U \quad V_1\right], &
	\Psi_{1}=\left[C \quad \underline{{V}}_{1}\right], &
	H_{1}=\left[
        		\begin{array}
	        [c]{cc}%
	        I & B_{1}\\
	        0 & {\underline{T}}_{1}%
	        \end{array}
	    \right], \\[10pt]
	\tilde{\Phi}_{1} = \left[\tilde{U} \quad \tilde{V}_1\right], &
	\tilde{\Psi}_{1}=\left[\tilde{C} \quad \underline{\tilde{{V}}}_{1}\right], &
	\tilde{H}_{1}=\left[
	        \begin{array}
	        [c]{cc}%
	        I & \tilde{B}_{1}\\
	        0 & \tilde{{\underline{T}}}_{1}%
	        \end{array}
	    	\right],
	\end{array}
\end{align*}
where $\underline{{V}}_{1}$ and $\underline{\tilde{V}}_{1}$ denote
$[V_1 \, v_{s+1}]$ and $[\tilde{V}_1 \, \tilde{v}_{s+1}]$
respectively.

\subsection{Constructing Bi-orthogonal $C_j$, $\tilde{C}_j$ and $C$, $\tilde{C}$}\label{secOrthoC}
We need to compute the matrices $C_j$ and $\tilde{C}_j$ such that $C_j \perp_b \tilde{C}_j$
at the end of each cycle.
After solving the generalized eigenvalue problem (\ref{finalGenEigValue}), we set (as initial choice) $U_{j} = \Phi_j W_j$, $\tilde{U}_{j} = \tilde{\Phi}_j \tilde{W}_j$, $C_{j} = A U_{j}$, and $\tilde{C}_{j} = A^* \tilde{U}_{j}$, and we compute the SVD
\begin{equation}\label{equationSVD}
    \tilde{C}_{j}^*C_{j} = M_j \Sigma_j N_j^*,
\end{equation}
such that $\sigma_1 \ge \sigma_2 \ge \ldots \ge \sigma_k \ge 0$. Given some tolerance {\tt tol} $> 0$, we pick $p$ such that $\sigma_p \ge$ {\tt tol} $> \sigma_{p+1}$ (with both $p = k$ and $p=0$ possible), and redefine $M_j = \left[m_1, \ldots, m_p\right]$ and $N_j = \left[n_1, \ldots, n_p\right]$, where $m_i$ and $n_i$ are the left and right singular vectors corresponding to $\sigma_i$. Next, we redefine
\begin{align}\label{equationForUC}
	\begin{array}
	[l]{ll}
	U_j = \Phi_j W_j N_j = [U_{j-1} \ V_j]W_j N_j, & \tilde{U}_{j} = \tilde{\Phi}_j \tilde{W}_j M_j = [{\tilde U}_{j-1} \ {\tilde V}_j] \tilde{W}_j M_j, \\[5pt]
	C_j = A U_j = A[U_{j-1} \ V_j]W_j N_j, & \tilde{C}_{j} = A^*\tilde{U}_j = A^*[{\tilde U}_{j-1} \ {\tilde V}_j] \tilde{W}_j M_j.
	\end{array}
\end{align}
By construction $\tilde{C}_j^* C_j$ is diagonal with real, positive coefficients.\footnote{For $p=0$, no recycle space would be selected. This has never occurred in our experience. Indeed, discarding even one pair of vectors is rare.}

Analogous to the above, at the start of each linear system (after the first), we need to compute $C$ and $\tilde{C}$ such that $C \perp_b \tilde{C}$, cf. (\ref{biorth}), and the diagonal matrix $\mathcal{D}_c = \tilde{C}^* C$ has real, positive coefficients. Taking initially for $U$ and $\tilde{U}$ the final matrices $U_j$ and $\tilde{U}_j$ from the previous pair of linear systems, we compute $C = AU$, $\tilde{C} = A^*\tilde{U}$, and compute the SVD $\tilde{C}^*C = M\Sigma N^*$. After this we proceed as for the computation of $C_j$ and $\tilde{C}_j$.

\subsection{Efficiently Setting up the Generalized Eigenvalue Problem}\label{secSimplify}
The main cost of setting up the generalized eigenvalue problem
(\ref{finalGenEigValue}) is in computing the matrices
\begin{eqnarray*}
\tilde{\Psi}_{j}^*\Psi_{j}=\left[
\begin{array}
    [c]{c}%
    \tilde{C}^*\\
    \tilde{C}_{j-1}^*\\
    \tilde{\Upsilon}_{j}^*%
\end{array}
\right] \left[
\begin{array}
[c]{ccc}%
C & C_{j-1}\text{ } & \Upsilon_{j}%
\end{array}
\right]  & =&
\left[\begin{array}
    [c]{ccc}%
    \mathcal{D}_c & \tilde{C}^*C_{j-1} & 0 \\
    \tilde{C}_{j-1}^* C & \Sigma_{j-1} &
    \tilde{C}_{j-1}^*\Upsilon_{j}\\
    0 & \tilde{\Upsilon}_{j}^* C_{j-1}\text{ } & I%
\end{array}\right],
\\
\tilde{\Psi}_{j}^* \Phi_{j}=\left[
\begin{array}
    [c]{c}%
    \tilde{C}^*\\
    \tilde{C}_{j-1}^*\\
    \tilde{\Upsilon}_{j}^*%
    \end{array}
\right] \left[
\begin{array}
[c]{cc}%
U_{j-1} & V_{j}%
\end{array}
\right]  & = &
\left[\begin{array}
    [c]{cc}%
    \tilde{C}^* U_{j-1} & 0\\
    \tilde{C}_{j-1}^* U_{j-1} & \tilde{C}_{j-1}^*V_{j}\\
    \tilde{\Upsilon}_{j}^*U_{j-1} & \overline{\underline{I}}
    \end{array}\right],
\end{eqnarray*}
where $\overline{\underline{I}}$ is the $s \times s$ identity matrix with an
extra row of zeros at the top and at the bottom.
Most of the blocks in these matrices can be constructed efficiently by exploiting
recurrences and various algebraic relations.

The bi-orthogonality
condition (\ref{biorth}) and the construction of
$C_j$ and $\tilde{C}_j$ (\ref{equationSVD})-(\ref{equationForUC})
give the following orthogonality relations.
\begin{equation}\label{biorthj}
	\tilde{C}_{j-2} \perp \Upsilon_{j}, \quad C_{j-2} \perp \tilde{\Upsilon}_{j}.
\end{equation}
Next, going from top-to-bottom and left-to-right,
we analyze each block of $\tilde{\Psi}_{j}^* \Psi_{j}$ and $\tilde{\Psi}_{j}^*\Phi_{j}$
in terms of its defining recurrences and simplify it using
(\ref{biorth}), (\ref{augBiLanExtend}), (\ref{equationSVD}), (\ref{equationForUC}), and (\ref{biorthj}).
Blocks whose efficient computation is obvious or has already been detailed are skipped,
and we focus on computations that are at least $O(n)$.
\begin{flalign*}
   \bullet \quad \tilde{C}^* C_{j-1} & = \tilde{C}^*AU_{j-1} = \left[
        \begin{array}
        [c]{cc}%
        \tilde{C}^*AU_{j-2} & \tilde{C}^*AV_{j-1}%
        \end{array}
        \right]W_{j-1}N_{j-1} & \\
    & = \left[
        \begin{array}
        [c]{cc}%
        \tilde{C}^*C_{j-2} & \tilde{C}^*\left(  C\hat{C}^{\ast
        }AV_{j-1}+\Upsilon_{j-1}\Gamma_{j-1}\right)
        \end{array}
        \right]W_{j-1}N_{j-1} & \\
    & = \left[
            \begin{array}
            [c]{cc}%
            \tilde{C}^*C_{j-2} & \mathcal{D}_cB_{j-1}
            \end{array}
            \right]W_{j-1}N_{j-1}, &
\end{flalign*}
For this first block, we describe its efficient computation in some detail.
The cost of computing $\tilde{C}^* C_{j-1}$ by direct multiplication is $O(k^2 n)$; so, it would be expensive.
However, the submatrix $\tilde{C}^*C_{j-2}$ is available from the previous cycle, and
$\mathcal{D}_c$ is a diagonal matrix independent of the cycle (so both must be computed at most
once per linear system). Furthermore, $B_{j-1}$ has been computed during the
(augmented) bi-Lanczos iteration.
Finally, the matrix-matrix product $[\tilde{C}^*C_{j-2} \  \ \mathcal{D}_cB_{j-1}]W_{j-1}N_{j-1}$ does
not involve any $O(n)$ operation. Hence, this block can be computed quite cheaply.
We give a brief overview of the cost of the RBiCG algorithm in Section \ref{resultsSec};
for a more detailed derivation, see \cite{Ahuja11_phd}.
\begin{flalign*}
    \bullet \quad \tilde{C}_{j-1}^* C
    & = M_{j-1}^* \tilde{W}_{j-1}^* \left[
            \begin{array}
            [c]{c}%
            \tilde{C}_{j-2}^*C\\
            {\tilde B}_{j-1}^* \mathcal{D}_c
            \end{array}
            \right]. &
\end{flalign*}
The derivation of this block is similar to the previous block.
As above, $\tilde{C}_{j-2}^*C$ is available from the previous cycle,
and $\tilde{B}_{j-1}$ has been computed during the bi-Lanczos iteration.
\begin{flalign*}
    \bullet \quad \tilde{C}_{j-1}^*\Upsilon_{j} &
    = \tilde{U}_{j-1}^*A\Upsilon_{j}
    = M_{j-1}^* \tilde{W}_{j-1}^* \left[
        \begin{array}
        [c]{c}%
        \tilde{U}_{j-2}^*\\
        \tilde{V}_{j-1}^*%
        \end{array}
        \right]  A\Upsilon_{j}
     = M_{j-1}^* \tilde{W}_{j-1}^* \left[
        \begin{array}
        [c]{c}%
        0 \\
        \tilde{V}_{j-1}^*A\Upsilon_{j}%
        \end{array}
        \right] & \\
    & = M_{j-1}^* \tilde{W}_{j-1}^* \left[
        \begin{array}
        [c]{c}%
        0 \\
        \left[\tilde{C}\check{C}^*A^*\tilde{V}_{j-1} + \tilde{\Upsilon}_{j-1}\tilde{\Gamma}_{j-1}\right]^*\Upsilon_{j}
        \end{array}
        \right] & \\
    & = M_{j-1}^* \tilde{W}_{j-1}^* \left[
        \begin{array}
        [c]{c}%
        0 \\
        \tilde{\Gamma}_{j-1}^*\tilde{\Upsilon}_{j-1}^*\Upsilon_{j}
        \end{array}
        \right], &
\end{flalign*}
where
\begin{align*}
    \tilde{\Gamma}_{j-1}^*\tilde{\Upsilon}_{j-1}^*
    \Upsilon_{j}
    & = \tilde{\Gamma}_{j-1}^*
    \left[
    \begin{array}
        [c]{c}%
        \tilde{v}_{(j-2)s}^*\\
        \tilde{v}_{(j-2)s+1}^*\\
        \vdots\\
        \tilde{v}_{(j-1)s}^*\\
        \tilde{v}_{(j-1)s+1}^*%
    \end{array}
    \right]
    \left[
    \begin{array}
        [c]{ccccc}%
        {v}_{(j-1)s} &
        {v}_{(j-1)s+1} &
        \cdots &
        {v}_{js} &
        {v}_{js+1}
    \end{array}
    \right] \\
    & = \tilde{\Gamma}_{j-1}^*
    \left[
    \begin{array}
        [c]{ccccc}%
        0 & 0 & \cdots & 0 & 0 \\
        \vdots & \vdots & \vdots & \vdots & \vdots \\
        0 & 0 & \cdots & 0 & 0 \\
        1 & 0 & \cdots & 0 & 0 \\
        0 & 1 & \cdots & 0 & 0
    \end{array}
    \right].
\end{align*}
\begin{flalign*}
    \bullet \quad \tilde{\Upsilon}_{j}^*C_{j-1}
    & = \left[
        \begin{array}
        [c]{cc}%
        0 &
        \tilde{\Upsilon}_{j}^*\Upsilon_{j-1}\Gamma_{j-1}
        \end{array}
        \right] W_{j-1} N_{j-1}, \quad\mbox{ where } & &
\end{flalign*}
\begin{align*}
    \tilde{\Upsilon}_{j}^*\Upsilon_{j-1}\Gamma_{j-1}
    & =
    \left[
    \begin{array}
        [c]{ccccc}%
        0 & 0 & \cdots & 1 & 0 \\
        0 & 0 & \cdots & 0 & 1 \\
        0 & 0 & \cdots & 0 & 0 \\
        \vdots & \vdots & \vdots & \vdots & \vdots \\
        0 & 0 & \cdots & 0 & 0
    \end{array}
    \right]\Gamma_{j-1}. &
\end{align*}
The derivation of this block is similar to the previous block.
\begin{flalign*}
    \bullet \quad \tilde{C}^*U_{j-1}
    & = \tilde{C}^*\left[
        \begin{array}
        [c]{cc}%
        U_{j-2} & V_{j-1}%
        \end{array}
        \right] W_{j-1} N_{j-1}
    = \left[
        \begin{array}
        [c]{cc}%
        \tilde{C}^*U_{j-2} & 0%
        \end{array}
        \right] W_{j-1}N_{j-1}, &
\end{flalign*}
where $\tilde{C}^*U_{j-2}$ is available from the previous cycle (such a block must be computed
at most once per linear system).
\begin{flalign*}
    \bullet \quad \tilde{C}_{j-1}^*U_{j-1}
    & = M_{j-1}^* \tilde{W}_{j-1}^*\left[
        \begin{array}
        [c]{c}%
        \tilde{U}_{j-2}^*\\
        \tilde{V}_{j-1}^*%
        \end{array}
        \right]  A\left[
        \begin{array}
        [c]{cc}%
        U_{j-2} & V_{j-1}%
        \end{array}
        \right] W_{j-1} N_{j-1} & \\
    & = M_{j-1}^* \tilde{W}_{j-1}^*\left[
        \begin{array}
        [c]{cc}%
        \tilde{C}_{j-2}^*U_{j-2} & \tilde{C}_{j-2}^*V_{j-1}\\
        \tilde{V}_{j-1}^*C_{j-2} & \tilde{V}^*_{j-1}AV_{j-1}%
        \end{array}
        \right]W_{j-1} N_{j-1}, &
\end{flalign*}
where $\tilde{C}_{j-2}^*U_{j-2}$ and $\tilde{C}_{j-2}^*V_{j-1}$ are
submatrices of $\tilde{\Psi}_{j-1}^*\Phi_{j-1}$, $\tilde{V}_{j-1}^*C_{j-2}$ is a
submatrix of $\tilde{\Upsilon}_{j-1}^*C_{j-2}$ and is available
from $\tilde{\Psi}_{j-1}^*\Psi_{j-1}$, and
$\tilde{V}^*_{j-1}AV_{j-1} = \tilde{T}_{j-1}$.
\begin{flalign*}
    & \bullet \quad \text{$\tilde{C}_{j-1}^*V_{j}$ is a submatrix of
    $\tilde{C}_{j-1}^*\Upsilon_{j}$, and hence, is available from $\tilde{\Psi}_{j}^*\Psi_{j}$}. &
\end{flalign*}
Therefore, only $\tilde{\Upsilon}_j^*U_{j-1}$ needs to be computed explicitly.


\section{Model Reduction}\label{secModRed}
Consider a single-input/single-output (SISO) linear time-invariant (LTI) system represented as
\begin{align}\label{equation:orginalModel}
    \begin{array}
        [l]{ll}%
        G:
        \begin{cases}
            E\,{\dot x}(t) = Ax(t) + bv(t) \\
            y(t) = c^*x(t),
        \end{cases}
    \end{array}
    ~~~{\rm or}~~~ \quad G(s) = c^*(sE - A)^{-1}b
\end{align}
where $E, A \in \mathbb{R}^{n \times n}$ and $b, c \in \mathbb{R}^{n}$. The time-dependent functions $v(t)$, $y(t)$: $\mathbb{R} \rightarrow \mathbb{R}$ are the input and output of $G(s)$, respectively, and $x(t): \mathbb{R} \rightarrow \mathbb{R}^n$ is the associated state.
In (\ref{equation:orginalModel}), $G(s)$ is the transfer function of the system: Let ${V}(s)$ and $Y(s)$ denote the Laplace transforms of $v(t)$ and $y(t)$, respectively. Then, the transfer function $G(s)$ satisfies
$Y(s) = G(s) V(s).$ By a common abuse of notation, we denote both the underlying dynamical system and its transfer function with $G$. The dimension of the underlying state-space, $n$, is called the dimension of $G$. Systems of the form (\ref{equation:orginalModel})
with extremely large state-space dimension $n$ arise in many applications; see \cite{antoulas:book} and \cite{KorR05}
for a collection of such examples.  Simulations in such large scale settings lead to overwhelming demands on computational resources. This is the main motivation for model reduction. The goal is to produce a surrogate model of much smaller dimension which provides a high-fidelity approximation of the input-output behavior of the original model $G$.
 Let $r \ll n$ denote the order of the reduced-model.
 The reduced-model is represented, similar to (\ref{equation:orginalModel}), as
 \begin{align}\label{equation:reducedModel}
    \begin{array}
        [l]{ll}%
        G_r(s):
        \begin{cases}
           E_r\, {\dot x}_r(t) = A_r x_r(t) + b_r v(t) \\
            y_r(t) = c_r^* x_r(t),
        \end{cases}
    \end{array} ~~~{\rm or}~~~G_r(s) = c_r^*(sE_r - A_r)^{-1}b_r
\end{align}
where $E_r,A_r \in \mathbb{R}^{r \times r}$ and $b_r, c_r \in \mathbb{R}^{r}$.
In this setting, the common approach is to construct reduced order
models via a Petrov-Galerkin projection. This amounts to choosing
two $r$-dimensional  subspaces  $\mathcal{V}_r$ and   $\mathcal{W}_r$ and matrices
$V_r \in \mathbb{R}^{n \times r}$ and $W_r \in \mathbb{R}^{n \times r}$
such that $\mathcal{V}_r = \mathrm{Range}(V_r)$ and $\mathcal{W}_r = \mathrm{Range}(W_r)$.
Then, we approximate
the full-order state $x(t)$ as $x(t) \approx V_r x_r(t)$ and enforce
the Petrov-Galerkin condition,
\begin{eqnarray*} \label{eqn:pg}
W_{\! r}^{*}\left(E V_{\! r}\dot{x}_r(t)-A V_{\! r}x_r(t) -b\,v(t)\right)={0},
\quad y_r(t) = c^* V_r x_r(t),
\end{eqnarray*}
 leading to a reduced-order model as in (\ref{equation:reducedModel}) with
\begin{equation}  \label{red_projection}
E_r = W_r^*  E V_r,~~~ 
A_r = W_r^* A V_r,~~~
b_r = W_r^* b,~~~{\rm and}~~~
c_r = V_r^*c.
\end{equation}
As (\ref{red_projection}) illustrates, the quality of the reduced model depends
solely on the selection of the two subspaces
$\mathcal{V}_r$ and   $\mathcal{W}_r$. In this paper, we will choose
$\mathcal{V}_r$ and $\mathcal{W}_r$ to enforce interpolation.
For other selections of $\mathcal{V}_r$ and   $\mathcal{W}_r$, we refer the reader to~\cite{antoulas:book}.

\subsection{Interpolatory Model Reduction}\label{secKryModRed}
For a given full-order model $G(s)$, the goal of interpolatory model reduction
is to construct  a reduced-order model $G_r(s)$ via  rational interpolation.
Here, we focus on Hermite interpolation. Given the full-order model
(\ref{equation:orginalModel}) and a collection of interpolation points
(also called shifts) $\sigma_i \in \mathbb{C}$, for $i=1,
\ldots,r$, we must construct a reduced-order system
by projection as in (\ref{red_projection})
such that $G_r(s)$ interpolates $G(s)$ and
its first derivative at selected interpolation points, i.e.,
\begin{equation*} \label{eqn:hermite}
G( \sigma_i) = G_r( \sigma_i) \qquad {\rm and} \qquad G'( \sigma_i) = G'_r( \sigma_i)~~~{\rm for}~~~i=1,\ldots,r.
\end{equation*}
Rational interpolation by projection was first proposed in \cite{Skeltonlate,Skelt1,Skelt2}.
How to obtain the required projection was derived in \cite{Grim} using the rational Krylov method \cite{ruhe1984rational}.
For the special case of Hermite rational interpolation,
the solution of the interpolatory model reduction problem is given in Theorem \ref{thm:rationalkrylov}.
For the more general case, we refer the reader to  \cite{Grim} and the recent survey \cite{Ant2010imr}.
\begin{theorem}\label{thm:rationalkrylov}
Given $G(s)= c^*(s E- A)^{-1}b$ and $r$ distinct points $ \sigma_1, \dots,  \sigma_{r} \in \mathbb{C}$, let
\begin{equation} \label{eqn:VrWr}
{V}_r=\lbrack ( \sigma_1 {E} - {A})^{-1} {b} \dots  ( \sigma_r {E} - {A})^{-1} {b}\rbrack ,
\qquad
 {W}_r^*=
\begin{bmatrix}
 {c}^*( \sigma_{1} {E}- {A})^{-1}\\
\vdots \\
 {c}^*( \sigma_{r} {E}- {A})^{-1}
\end{bmatrix}.
\end{equation}
Using (\ref{red_projection}), define the reduced-order model $G_r(s) = c_r^*(s E_r - A_r)^{-1}b_r$.
Then $G(\sigma_i)=G_r(\sigma_i)$ and $G'(\sigma_i)=G'_r(\sigma_i)$, for $i=1,\dots, r$,
provided that
$\sigma_i E - A$ and $\sigma_i E_r - A_r$ are invertible for $i=1,\ldots,r$.
\end{theorem}

Theorem \ref{thm:rationalkrylov} shows how to solve the
interpolatory model reduction problem via projection
for given shifts. However, it does not provide a strategy for choosing good/ optimal interpolation points.
Recently, this issue has been resolved for the special case of optimality
in the $\mathcal{H}_2$ norm \cite{gugercin:base}.
The $\mathcal{H}_2$ norm of the dynamical system
$G(s)$  is defined as
\begin{equation*}
	\left \| G \right \|_{\Htwo} = \left(\frac{1}{2\pi} \int_{-\infty}^\infty \mid G(\jmath a) \mid^2 d a \right)^{1/2}.
\end{equation*}
The $\Htwo$ norm of $G$ is the $2-\infty$ induced norm of the underlying convolution operator.
Then, for any $v\in L^2(\mathbb{R}^+)$,
$\| y- y_{r} \|_{L^\infty} \leq  \| G- G_{r} \|_{\Htwo} \| v \|_{L^2}.$
To ensure that the output error $y- y_{r}$ is small in $ L^\infty(\mathbb{R}^+)$  uniformly over
all inputs $v$, say, with $\|v\|_{L^2}\leq 1$, we seek a reduced system $G_r$ that
makes  $ \| G - G_{r} \|_{\Htwo}$ small.
This leads to the \emph{optimal $\Htwo$ model reduction problem}:  Given  $G(s)$,
and a  reduced order $r < n$,
find  ${\displaystyle G_r(s)}$
that solves
\begin{equation} \label{optimalHtwo}
	\| G - G_r \|_{\Htwo} = \min \limits_{dim(\hat{G}_r) = r}  \left\| G-\hat{G}_r \right\|_{\Htwo}.
\end{equation}
This problem has been studied extensively
\cite{meieriii1967approximation,wilson1970optimum, gugercin:base,spanos1992anewalgorithm,vandooren2008hom,gugercin2005irk,beattie2007kbm,beattie2009trm}.
It is a non-convex optimization problem,
which makes finding the global minimum, at best, a hard task. Hence, the common approach is to construct reduced-order models that satisfy, for an interpolatory model reduction framework,
the following first-order necessary conditions.
\begin{theorem}(\cite{meieriii1967approximation, gugercin:base}) \label{thm:h2optimal}
Given $G(s)$, let $G_r(s) = c_r^*(s E_r - A_r)^{-1}b_r$ be an $\Htwo$-optimal
reduced order model of order $r$, and let
$\hat{\lambda}_1,\dots,\hat{\lambda}_r$ denote the poles of $G(s)$.
Then
\begin{equation} \label{eqn:h2cond}
G(-\hat{\lambda}_i)=G_r(-\hat{\lambda}_i)~~~{\rm and}~~~
G'(-\hat{\lambda}_i)=G_r'(-\hat{\lambda}_i)~~~{\rm for}~~~
 i=1,\dots,r.
\end{equation}
\end{theorem}
So, the $\Htwo$ optimal approximant $G_r(s)$ is a Hermite interpolant to $G(s)$
at the mirror image of its poles.
These poles, the optimal interpolation points, are not known \emph{a priori}.
Hence, the iterative rational Krylov algorithm ({IRKA}) \cite{gugercin:base},
starting from an initial selection of interpolation  points, iteratively corrects
the interpolation points until (\ref{eqn:h2cond}) is satisfied.
Algorithm 3 outlines IRKA; for details, see \cite{gugercin:base}.
\begin{figure}[t] {\bf Algorithm 3.} {\it IRKA (~\cite{gugercin:base})} \\
        1. Make an initial shift selection $\sigma_i$ for $i = 1, \ldots, r$, \\
        2. $V_r = [(\sigma_1 E - A)^{-1} b,\ \ldots,\ (\sigma_r E - A)^{-1} b]$ , \\
        3. $W_r = [(\sigma_1 E- A)^{-*} c,\ \ldots,\ (\sigma_r E - A)^{-*} c]$, \\
        4. while (not converged) \\
            $\diamond$ \; \; $A_r =W_r^*AV_r$, $E_r =W_r^*AV_r$,  \\
            $\diamond$ \; \; $\sigma_i \leftarrow -\lambda_i(A_r,E_r)$ for $i = 1, \ldots, r$, \\
            $\diamond$ \; \; $V_r = [(\sigma_1 E - A)^{-1} b,\ \ldots,\ (\sigma_r E - A)^{-1} b]$ , \\
            $\diamond$ \; \; $W_r = [(\sigma_1 E - A)^{-*} c,\ \ldots,\ (\sigma_r E - A)^{-*} c]$, \\
        5. $A_r =W_r^*AV_r$, $E_r =W_r^*AV_r$, $b_r =W_r^*b$, $c_r = V_r^*c_r$.
\end{figure}

\subsection{Advantages of Approximating Solutions using a
Petrov-Galerkin Framework in Interpolatory Model Reduction} \label{sec:pg}

The main cost in IRKA is solving multiple linear systems
to compute $V_r$ and $W_r$. If the dimension of the state-space, $n$, is large, these systems
are generally solved only approximately by an iterative solver. In this context, it is
important to asses the accuracy of the computed reduced order model; that is,
given the shifts, how accurately the Hermite interpolation problem is solved.
This question was
studied extensively in \cite{beattie2010isi}. One of the major results,
outlined below for our particular case,
provides the main motivation for solving the linear systems associated with
the corresponding columns of $V_r$ and $W_r$ as pairs of dual linear systems
(in the terminology of Section \ref{introSec}) using BiCG or
RBiCG.

Let $\hat{v}_j$ and $\hat{w}_j$, for $j=1,\ldots,r$, denote the approximate solutions
of $(\sigma_j E -A) v_j = b$ and  $(\sigma_j E -A)^* w_j = c$, respectively,
with residuals $\eta_j = (\sigma_j E -A) \hat{v}_j - b$ and
$\xi_j = (\sigma_j E -A)^* \hat{w}_j - c$.
Furthermore, let $\hat{v}_j$, $\hat{w}_j$, $\eta_j$, and $\xi_j$ satisfy
the Petrov-Galerkin condition that there exist spaces $\cal P$ and $\cal Q$ such that
$\hat{v}_j \in {\cal P}$, $\hat{w}_j \in {\cal Q}$, $\eta_j \perp {\cal Q}$,
and $\xi_j \perp {\cal P}$.  Define the approximate solution matrices
($\hat{V}_r$ and $\hat{W}_r$), the residual matrices ($R_b$ and $R_c$),
and the rank-2r matrix ($F_{2r}$) as follows:
\begin{align*}
	\begin{array}
	[l]{ll}
	\hat{V}_r =  [\hat{v}_1 \, \hat{v}_2 \, \ldots \, \hat{v}_r] , &
    	\hat{W}_r =  [\hat{w}_1 \, \hat{w}_2 \, \ldots \, \hat{w}_r], \\
	R_b =  [\eta_1 \, \eta_2 \, \ldots \, \eta_r], &
    	R_c =  [\xi_1 \, \xi_2 \, \ldots \, \xi_r] , \\[5pt]
    	F_{2r} = R_b (\hat{W}_r^*\hat{V}_r)^{-1}\hat{W}_r^* + \hat{V}_r^*(\hat{W}_r^*\hat{V}_r)^{-1}R_c. &
	\end{array}
\end{align*}
Also, define the inexact reduced-order order quantities
\begin{equation*}  \label{red_projection_inexact}
	\hat{A}_r = \hat{W}_r^* A \hat{V}_r,~~~
	\hat{E}_r = \hat{W}_r^*  E \hat{V}_r,~~~ 
	\hat{b}_r = \hat{W}_r^* b,~~~{\rm and}~~~
	\hat{c}_r = \hat{V}_r^*c.
\end{equation*}
Then, the computed reduced-order model $\hat{G}_r(s) = \hat{c}_r^*(s \hat{E}_r-\hat{A}_r)^{-1}\hat{b}_r$
{\it exactly} interpolates the perturbed
full-order model $\hat{G}(s) = c^* (sE-(A+F_{2r}))^{-1}b$, i.e.,
$$\hat{G}(\sigma_i)=\hat{G}_r(\sigma_i)~~~{\rm and}~~~\hat{G}'(\sigma_i)=\hat{G}'_r(\sigma_i),~~~{\rm for}~~~i=1,\dots, r. $$
Hence, iteratively solving the linear systems while satisfying the Petrov-Galerkin
condition above yields a backward error for the interpolatory model reduction that is
bounded by $\|F_{2r}\|$, which is governed by the
norms of the residuals. The latter are easily controlled in the iterative solver.
For details, we refer the reader to \cite{beattie2010isi}.

The easiest way to satisfy the Petrov-Galerkin condition above is by solving the
dual pairs of linear systems using BiCG.
Hence, BiCG is particularly suitable for solving the linear systems in
IRKA (even for symmetric positive definite matrices).
However, as IRKA leads to a sequence of dual linear systems, the RBiCG algorithm can
be used to reduce the total run time for solving all linear systems.
Moreover, if we solve the dual pairs of linear systems arising in IRKA by RBiCG,
the Petrov-Galerkin condition is still satisfied.
Hence, the resulting reduced-order model will be an optimal
$\Htwo$ approximation to a nearby full-order model.

\subsection{IRKA using RBiCG}\label{secIRKAusingRBiCG}
IRKA usually converges rather fast \cite{gugercin:base}. Hence, after one or
a few initial steps, the interpolations points from one step of IRKA to the next
do not change substantially. Moreover, for many cases, the change of the
(appropriately ordered)  \{$\sigma_i$\} from one column of $V_r$ (and $W_r$) to
the next is also modest. Therefore, IRKA is expected to gain significantly from recycling.

For the special case of $E=I$ in (\ref{equation:orginalModel}), alternative solution
approaches might be advantageous, as one can solve the linear systems for
multiple shifts at once~\cite{freund1993shifted, frommer:bicgstabl,
jegerlehner:shifted, vandenEshof:shifted}. Combining these strategies with a
Petrov-Galerkin framework does not seem complicated. Effective strategies for
Krylov subspace recycling for solving systems of the type,
$(\sigma_i I - A)v_i = b$, for multiple shifts at once, as well as for multiple right
hand sides, was discussed in~\cite{kilmer2006dot}. For most model reduction
applications, however, $E \neq I$.

There are three strategies for recycling Krylov subspaces in IRKA.
For the first strategy, consider two consecutive steps of IRKA, say step $m$
and $m+1$
(iterations $m$ and $m+1$ of the while loop in Algorithm~3), with
shifts $\sigma^{(m)}_i$, for $i = 1, \ldots, r$ and
\begin{equation}\label{equation:pStage}
    \begin{array}
    [l]{l}
    V^{(m)}_r= [(\sigma^{(m)}_1 E - A)^{-1} b,\ \ldots,\ (\sigma^{(m)}_r E - A)^{-1} b], \\
    W^{(m)}_r = [(\sigma^{(m)}_1 E - A)^{-*} c,\ \ldots,\ (\sigma^{(m)}_r E - A)^{-*} c].
    \end{array}
\end{equation}
at step $m$ and with shifts $\sigma^{(m+1)}_i$, for $i = 1, \ldots, r$ and
\begin{equation*}\label{equation:pplusStage}
    \begin{array}
    [l]{l}
     V^{(m+1)}_r= [(\sigma^{(m+1)}_1 E - A)^{-1} b,\ \ldots,\ (\sigma^{(m+1)}_r E - A)^{-1} b], \\
    W^{(m+1)}_r = [(\sigma^{(m+1)}_1 E - A)^{-*} c,\ \ldots,\ (\sigma^{(m+1)}_r E - A)^{-*} c]
    \end{array}
\end{equation*}
at step $m+1$ of IRKA.
One can recycle Krylov subspaces from the $i^{th}$ column of $V_r^{(m)}$ and
$W_r^{(m)}$ to the $i^{th}$ column of $V_r^{(m+1)}$ and $W_r^{(m+1)}$.
That is, from solving the pair of linear systems
$$(\sigma^{(m)}_i E - A)v^{(m)}_i = b, \quad (\sigma^{(m)}_i E - A)^{*}w^{(m)}_i =c,$$
to solving the pair of linear systems
$$(\sigma^{(m+1)}_i E - A)v^{(m+1)}_i = b, \quad (\sigma^{(m+1)}_i E - A)^{*}w^{(m+1)}_i =c,$$
where $i = 1, 2, \ldots, r$. This strategy for recycling strategy is useful
when the change in a shift from one IRKA step to the next is small.

For the second strategy, consider a single IRKA step.
One can recycle selected Krylov subspaces from solving for one pair of
columns of $V_r$ and $W_r$ to the next pair of columns across all the columns of the
matrices $V_r$ and $W_r$. In the third strategy, the first two recycling strategies are
combined. We describe one such combination. Consider solving the
system $(\sigma^{(m+1)}_i E - A)v^{(m+1)}_i = b$ and its dual system. >From a set of
previously generated recycle spaces (distinguished by their shifts), one can
pick the recycle space from the system with the smallest relative change in $\sigma$
(and less than a relative tolerance). This would ensure that the linear system from which
the recycle space has been generated is close to the current one.
A natural pool from which to
pick the $\sigma$ defining the recycle space would be
$\sigma^{(m)}_1, \ldots, \sigma^{(m)}_r, \sigma^{(m+1)}_1, \ldots, \sigma^{(m+1)}_{i-1}$.
The second and third recycling strategies are useful when the shifts at an
IRKA step are clustered.

For the experiments in this paper, $r$ is small, and so the shifts at any particular
IRKA step are spread far apart. Hence, we follow the first strategy. That is, for
every shift, we recycle Krylov subspaces from one IRKA step to the next. In general,
the linear systems corresponding to the relatively large shifts converge fast, and
so recycling Krylov subspaces is not useful for them. Therefore, we carry out
recycling only for selected, small shifts. We give more details in
Section \ref{secModRedResults}.

\subsection{Previous Work in Recycling for Model Reduction}\label{secModRedRecycle}
Recycling for interpolatory model reduction in the Galerkin setting, i.e., with $W_r = V_r$,
has been considered in \cite{benner:recycle} and \cite{feng:proceedings}.
In this setting, there are no dual systems to solve,
and therefore approaches based on GCR~\cite{eisenstat:gcr}
and GMRES~\cite{saad:gmres} are
considered, respectively, for a sequence of (single) linear systems,
as opposed to our approach based on BiCG for a sequence of dual linear systems.
Also in other respects, the approach for
improving the linear solver and the model reduction context are quite different from here.
In~\cite{benner:recycle}, the focus is on efficiently solving linear systems
with a fixed coefficient matrix and
multiple right hand sides ($Ax^{(j)} = b^{(j)}$),
recycling descent vectors (in GCR). Furthermore, the authors target model reduction with a
single interpolation point but interpolating higher derivatives.

\section{Results}\label{resultsSec}
We first give a brief overview of the overhead in RBiCG.
We focus on components with at least $O(n)$ cost, where $n$ is the dimension of the linear system.
Furthermore, $k$ is the number of basis vectors in the primal (or dual) recycle space,
and $s$ is the number of iterations in a cycle.
For every iteration, there is an extra cost of $(8k+2)n$ flops, mostly from orthogonalizations.
At the end of each cycle, there is an extra cost of $(14k^2 + 6ks +16k+4)n$ flops,
mostly from setting up the generalized eigenvalue problem and computing biorthogonal
$C_j$ and $\tilde{C}_j$.
Once per linear system, there is an extra cost of
$(10k^2 + 28k + 14)n$ flops, mostly from computing biorthogonal $C$ and $\tilde{C}$.
A more detailed discussion of the
overhead is given in \cite{Ahuja11_phd}.
Note that $s$ and $k$ are much smaller than $n$. For recycling to be beneficial,
the savings in iterations should be sufficient to make up for the overhead.
Further in this section, we show that the reduction in the number of
iterations for (a pair of) linear systems
may be as high as $70$\%. For our model reduction test problem,
we show that computing a reduced model
without recycling takes about 50\% more time than with recycling.

We test RBiCG on a convection-diffusion problem and on IRKA for interpolatory
model reduction. All experiments are done using Matlab.

\subsection{Convection-Diffusion}\label{secConvDiffResults}
\begin{figure}
\begin{tabular}{cc}
\includegraphics[width=2.5in]{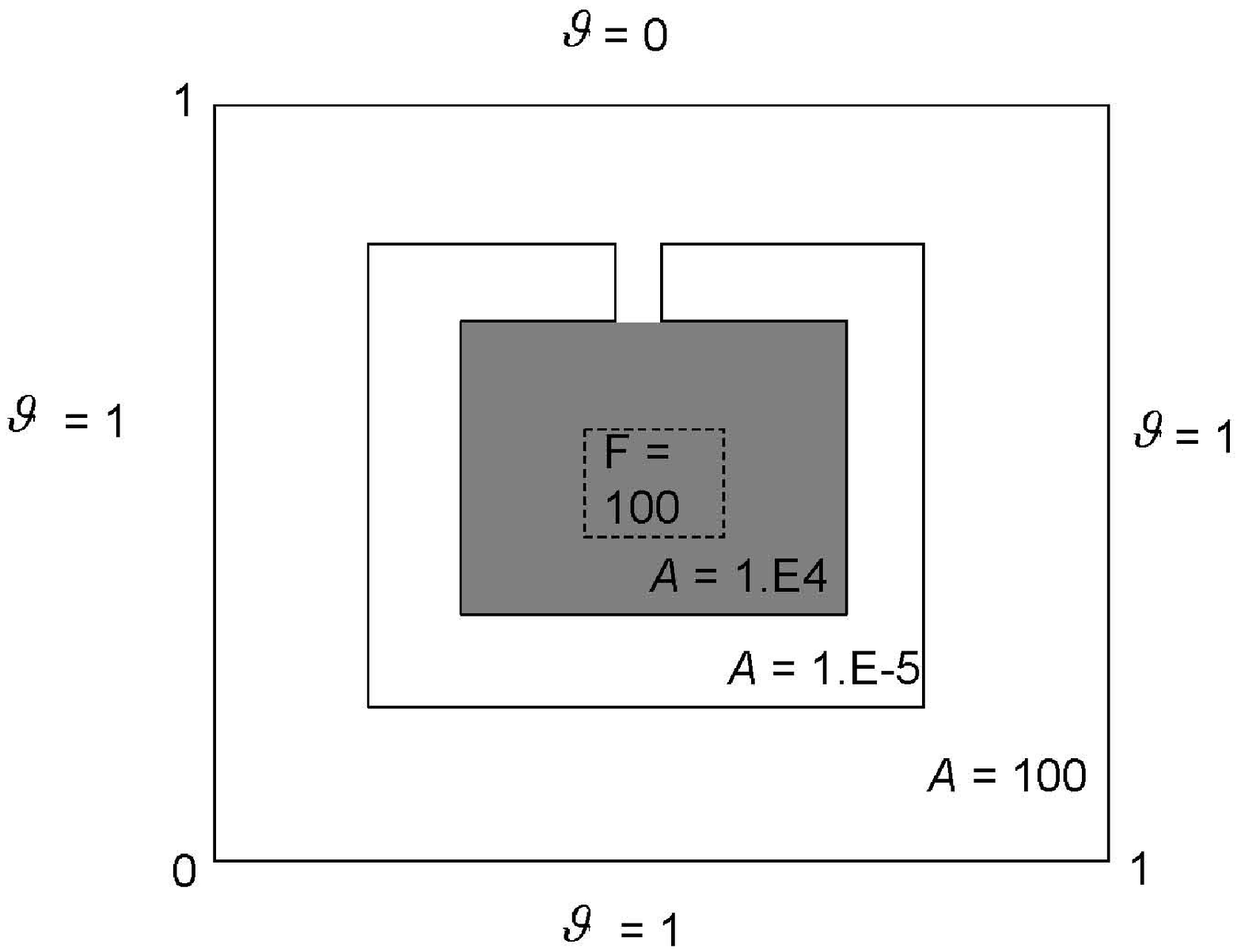} &
\includegraphics[width=2.5in]{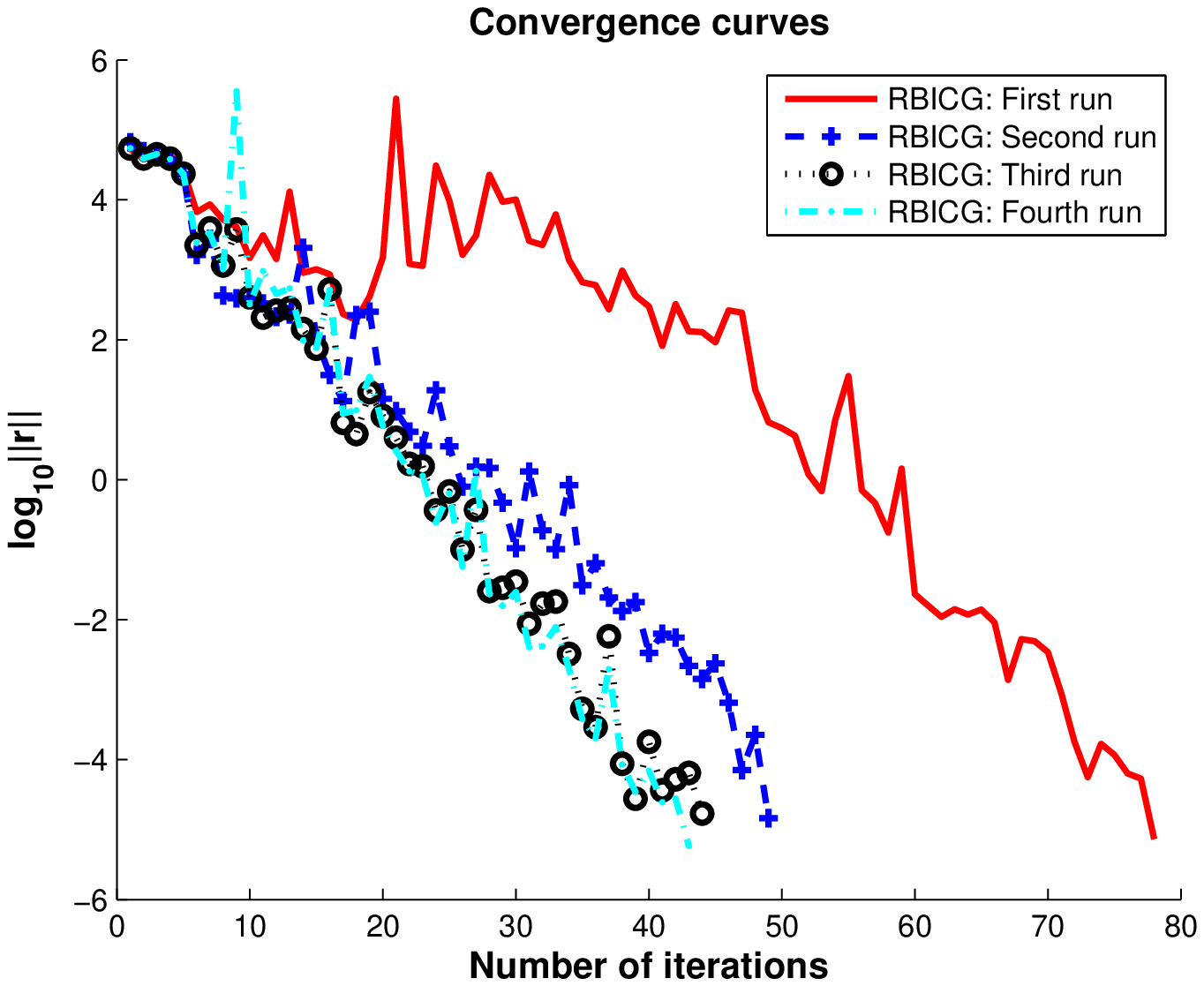}
\end{tabular}
    \caption{RBiCG for a convection-diffusion problem. (a) The coefficients of the PDE. (b) Convergence
    for preconditioned RBiCG with $s=40$ and $k=10$ for the primary system solved four times to analyze
    convergence improvement as the recycle space improves.}
    \label{figure:conv-diff}
\end{figure}

To analyze RBiCG, we use the linear system obtained by finite difference
discretization of the partial differential equation
\begin{align*}
    -(\mathcal{A}\vartheta_x)_x - (\mathcal{A}\vartheta_y)_y + \mathcal{B}(x,y)\vartheta_x = \mathcal{F},
\end{align*}
with $\mathcal{A}$ as shown in Figure \ref{figure:conv-diff} (a), $\mathcal{B}(x,y) = 2e^{2(x^2 + y^2)}$, and $\mathcal{F} = 0$ everywhere except in a small square in the center where $\mathcal{F} = 100$~\cite{vorst:bicgstab};
see Figure~\ref{figure:conv-diff}(a). The domain is $(0, 1) \times (0, 1)$ with Dirichlet boundary conditions
\begin{align*}
    \begin{array}
    [l]{l}
    \vartheta(0, y) = \vartheta(1, y) = \vartheta(x, 0) = 1, \mbox{ and }
    \vartheta(x, 1) = 0.
    \end{array}
\end{align*}
We use the second order central difference scheme with a mesh width of $h = 1/64$, giving
a nonsymmetric linear system of $3969$ unknowns.
The convergence is similar for a problem that is four times larger.
To enable further analysis, we give results for this smaller system size.
The primary system right-hand side comes from the PDE. We take the vector of all zeros
as the dual system right-hand side.
In this case, we are concerned only about the primary system.

To analyze RBiCG, we solve the (same) dual linear systems four times.
The recycle space generated during the first run is used for solving the same dual systems a second time, further improving the recycle space, and so on.
This is a useful approach for analyzing how well Krylov subspace recycling works,
as it excludes the effects of changing matrices and of right-hand sides having different expansions in the eigenvector basis~\cite{mike:recycle}.
Hence, it provides an indication for reasonable sequences of systems
how fast the recycle spaces converge and how much
recycling approximate invariant subspaces is likely to improve convergence.
For this experiment, we take $s = 40$ and $k = 10$. These are chosen based on experience with other recycling algorithms~\cite{mike:recycle, shun:base}.
The relative tolerance for RBiCG is taken as $10^{-8}$. The initial guess (for both systems)
is a vector of all ones. The linear systems are split-preconditioned by a Crout version of the ILUT preconditioner with a
drop tolerance of $0.05$~\cite{saad:book}. The generated recycle spaces pertain to the preconditioned linear systems.

Figure \ref{figure:conv-diff}(b) shows the convergence improvement of RBiCG,
as it solves the primary system multiple times. For the second run, the reduction
in iterations is around 35\%. The convergence improves further with each run. Next, we present a brief analysis of the generated recycle spaces. In Table~\ref{table:eigenvectorsConvDiff}, we give the cosines of the principal angles between
the primary (dual) recycle space and the right (left) invariant subspace associated with the eight
eigenvalues of smallest magnitude. As for the recycle spaces, the invariant subspaces are computed for the preconditioned operator.
As the recycle space improves, the principal angles to tend to zero, and so the cosines tend to one. The table shows that with only a few runs, RBiCG accurately approximates increasingly larger
subspaces of the invariant subspace. As a result, we see faster convergence for every new run.
%

\begin{table}[bht]
   \footnotesize
    \begin{center}
        \setlength{\extrarowheight}{1.5pt}
        \begin{tabular}{|m{0.5in}m{0.5in}m{0.5in}|m{0.5in}m{0.5in}m{0.5in}|}
         	   \hline
            \multicolumn{3}{|c|}{Primary System} & \multicolumn{3}{c|}{Dual System} \\
         	   \hline
	   Start of & Start of & Start of & Start of & Start of & Start of \\
            Run 2 & Run 3 & Run 4 &  Run 2 & Run 3 & Run 4 \\
         	   \hline
	   1.0000	&	1.0000      &	1.0000 	& 	1.0000		& 	1.0000		& 	1.0000		\\
	   1.0000	&	1.0000      &	1.0000 	& 	1.0000		& 	1.0000		& 	1.0000		\\
	   1.0000	&	1.0000      &	1.0000 	& 	1.0000		& 	1.0000		& 	1.0000		\\
	   0.9896	&	1.0000	&	1.0000 	& 	0.9950		& 	1.0000		& 	1.0000		\\
	   0.3832	&	1.0000	&	1.0000 	& 	0.9884		& 	1.0000		& 	1.0000		\\
	   0.1452	&	0.9983	&	1.0000 	& 	0.7864		& 	0.9844		& 	1.0000		\\
	   0.0988	&	0.9437	&	0.9970 	& 	0.6070		& 	0.9206		& 	0.8141		\\
	   0.0300	&	0.1869	&	0.9567 	& 	0.4749		& 	0.4118		& 	0.4721		\\
         	   \hline
        \end{tabular}
    \end{center}
    \caption{Convergence of the recycle space for the convection-diffusion problem as measured by
    the cosines of the principal angles between the primary (dual) recycle space and the right (left)
    invariant subspace associated with the eight eigenvalues of smallest magnitude.}
    \label{table:eigenvectorsConvDiff}
\end{table}

\subsection{Model Reduction}\label{secModRedResults}
Our test dynamical system is a semi-discretized heat transfer problem for determining the optimal cooling of steel profiles~\cite{penzl1999mds, benner2004lcp, saak2004lqr}.
We will refer to this model
as the {\it rail model}~\cite{penzl1999mds}. The rail model  has seven inputs and six outputs. Since we focus on SISO systems in this paper, we choose a SISO subsystem corresponding to the second input and sixth output.
The rail model is available with $1357$, $5177$, $20209$, and $79841$ unknowns, depending on the mesh size.


As convergence tolerance for IRKA we use a
relative change in the shifts of less than  $10^{-6}$.
The matrices $A$ and $E$ of (\ref{equation:orginalModel}) are symmetric
negative definite and symmetric positive definite (SPD), respectively. Since our
shifts are real and positive at every IRKA step, $(\sigma_i^{(m)} E - A)$ is always SPD.
Nevertheless, RBiCG is advantageous here because of the backward error formulation
discussed in Section \ref{sec:pg}.
We carry out two sets of experiments that differ in the dimension, $r$, of the reduced models.
We also vary the frequency of computing a recycle space,
as a recycle space can be effective for multiple consecutive systems~\cite{mike:recycle, kilmer2006dot}
and updating it may be expensive.

\begin{figure}[tbh]
    \centering
    \includegraphics[scale=0.55]{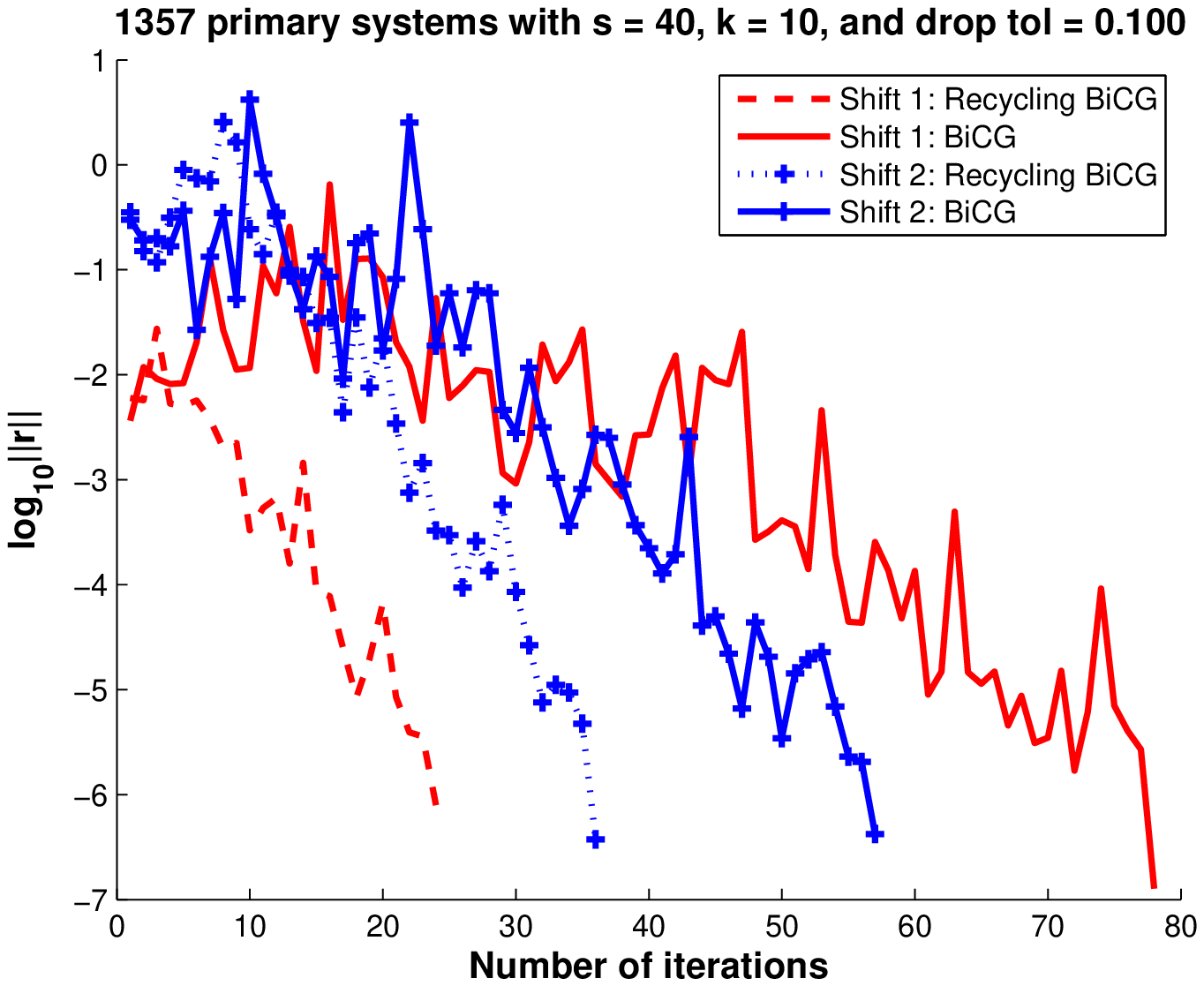}
    \caption{Convergence of preconditioned RBiCG at the $3^\text{rd}$ IRKA step for the $n=1357$ rail model, with
    $s = 40$, $k=10$, and the preconditioner drop tolerance is $0.1$. }
    \label{figure:irka1300}
\end{figure}

\begin{figure}[tbh]
    \centering
    \includegraphics[scale=0.55]{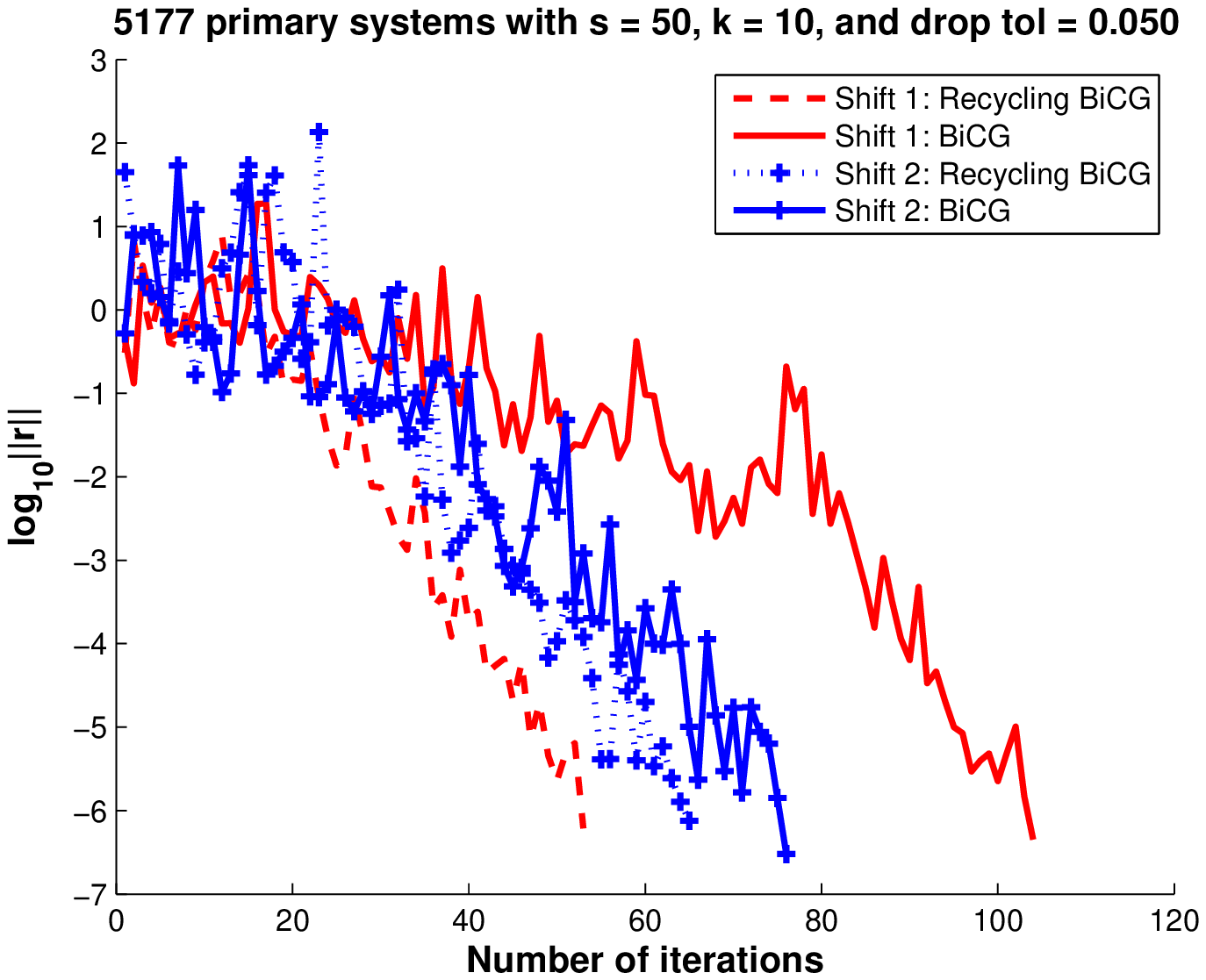}
    \caption{Convergence of preconditioned RBiCG at the $2^\text{nd}$ IRKA step for the $n=5177$ rail model, with $s=50$, $k=10$, and the preconditioner drop tolerance is $0.05$. }
    \label{figure:irka5k}
\end{figure}

\begin{figure}
    \centering
    \includegraphics[scale=0.55]{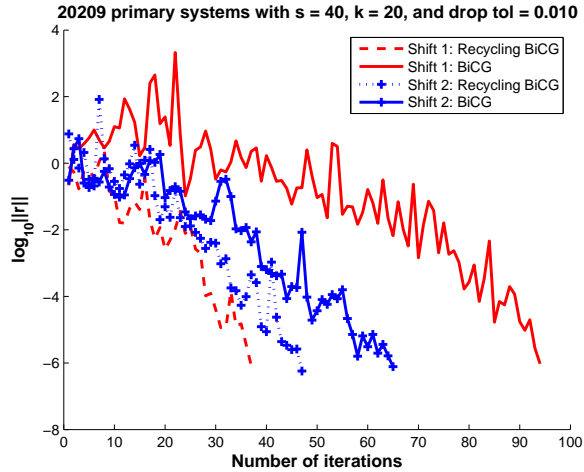}
    \caption{Convergence of preconditioned RBiCG at the $2^\text{nd}$ IRKA step for the $n=20209$ rail model, with $s=40$, $k=20$, and the preconditioner drop tolerance is $0.01$. }
    \label{figure:irka20K}
\end{figure}

\begin{figure}
    \centering
    \includegraphics[scale=0.55]{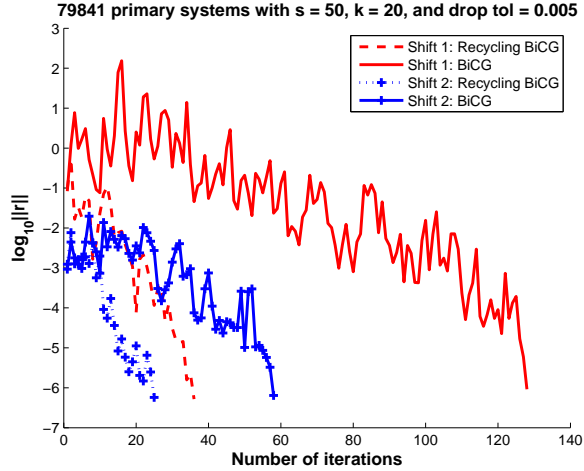}
    \caption{Convergence of preconditioned RBiCG at the $3^\text{rd}$ IRKA step for the $n=79841$ rail model,
    with $s=50$, $k=20$, and the preconditioner drop tolerance is $0.005$.}
    \label{figure:irka80k}
\end{figure}

We implement the first recycling strategy from Section \ref{secIRKAusingRBiCG} for a few selected shifts.
As for the convection-diffusion example, the recycling parameters $s$ and $k$ are chosen based
on experience with other recycling algorithms~\cite{mike:recycle, shun:base}. If a pair of
linear systems converges in fewer than $s$ iterations, the recycle space is not updated, and
we use the previous recycle space for the next pair of systems in the sequence.
The relative convergence tolerance for the iterative solves and the tolerance for
constructing $\tilde{C}_j$ and $C_j$ in Section \ref{secOrthoC} are
taken as $10^{-6}$. The linear systems are split-preconditioned with an
incomplete LU preconditioner with threshold and pivoting (ILUTP)~\cite{saad:book}.
The drop tolerance varies per problem to avoid ill-conditioning;
see Figures~\ref{figure:irka1300} -- \ref{figure:irka80k}.
The initial guess of the preconditioned system is the solution vector from the
previous preconditioned system in the sequence.
For the first IRKA step, we take a vector of all zeros as the initial guess.
In general, a better initial guess may be based on knowledge of the system and
aim to avoid orthogonal initial residuals (Algorithm 1 Step 2; Algorithm 2 Step 3).

For the first set of experiments, we reduce the models to $r = 6$ degrees of freedom,
with $1.00 \times 10^{-5}$, $1.38 \times 10^{-4}$, $1.91 \times 10^{-3}$, $2.63 \times 10^{-2}$,
$3.63 \times 10^{-1}$, and $5.01$ as initial shifts.

We compute a recycle space at every IRKA step. The results
for the primary systems at a particular IRKA step (given in the caption) are
given in Figures~\ref{figure:irka1300} -- \ref{figure:irka80k}.
The graphs for the other IRKA steps are similar, as are the graphs
for the dual systems. We carry out recycling for the smallest two shifts.
Each figure has two solid curves for the linear systems solved without recycling
and two dashed--dotted curves for those solved with recycling.
It is evident that recycling significantly reduces the number of iterations.
The savings in iterations are as high as 70\% per system.
As discussed in Section \ref{secIRKAusingRBiCG}, convergence for the
remaining four (larger) shifts is rapid, so recycling Krylov subspaces is not useful for these.

\begin{table}[tbh]
   \footnotesize
    \begin{center}
        \setlength{\extrarowheight}{1.5pt}
        \begin{tabular}{|m{0.5in}m{0.5in}|m{0.5in}m{0.5in}||m{0.5in}m{0.5in}|m{0.5in}m{0.5in}|}
         	   \hline
            \multicolumn{4}{|c||}{Primary System} & \multicolumn{4}{c|}{Dual System} \\
         	   \hline
            \multicolumn{2}{|c|}{IRKA Step 1} & \multicolumn{2}{c||}{IRKA Step 2} &  \multicolumn{2}{c|}{ IRKA Step 1} & \multicolumn{2}{c|}{IRKA Step 2} \\
            \multicolumn{2}{|c|}{$\sigma_1 = 1.000 \times 10^{-5}$} & \multicolumn{2}{c||}{$\sigma_1 = 1.834 \times 10^{-5}$} &  \multicolumn{2}{c|}{$\sigma_1 = 1.000 \times 10^{-5}$} & \multicolumn{2}{c|}{$\sigma_1 = 1.834 \times 10^{-5}$} \\         	   \hline
            End of 	& End of  	& Start of 	& End of 	& End of  	& End of  	& Start of 	& End of   \\
            Cycle 1	& Cycle 2	& Cycle 1 & Cycle 1  & Cycle 1 	& Cycle 2 	& Cycle 1 & Cycle 1  \\
         	   \hline
	   1.0000	&	1.0000      &	1.0000 	& 	1.0000		& 	1.0000		& 	1.0000	&	1.0000	&	1.0000	\\
	   1.0000	&	1.0000      &	1.0000 	& 	1.0000		& 	1.0000		& 	1.0000	&	1.0000	&	1.0000	\\
	   1.0000	&	1.0000      &	1.0000 	& 	1.0000		& 	0.9997		& 	1.0000	&	1.0000	&	1.0000	\\
	   0.9987	&	1.0000	&	1.0000 	& 	1.0000		& 	0.9765		& 	1.0000	&	1.0000	&	1.0000	\\
	   0.9321	&	1.0000	&	1.0000 	& 	1.0000		& 	0.4936		& 	1.0000	&	1.0000	&	1.0000	\\
	   0.2257	&	1.0000	&	0.9998 	& 	0.9999		& 	0.0844		& 	0.9995	&	0.9997	&	0.9998	\\
	   0.0260	&	0.9997	&	0.9996 	& 	0.9997		& 	0.0231		& 	0.9945	&	0.9945	&	0.9989	\\
	   0.0072	&	0.7813	&	0.7799 	& 	0.9932		& 	0.0068		& 	0.3439	&	0.3423	&	0.9876	\\
         	   \hline
        \end{tabular}
    \end{center}
    \caption{\label{table:eigenvectorsRail} Convergence of the recycle space for the sequence of linear systems
    corresponding to the $5177$ rail model and the smallest shift, as measured by
    the cosines of the principal angles between the primary (dual) recycle space and the right (left)
    invariant subspace associated with the eight eigenvalues of smallest magnitude.
    The third column 
    corresponds to the dashed convergence curve in Figure \ref{figure:irka5k}.}
\end{table}
Next, we analyze the recycle space generated during the first two IRKA steps for
the order $5177$ rail model corresponding to the smallest shift.
In Table \ref{table:eigenvectorsRail}, we give the cosines of principal angles between
the recycle space and the invariant subspace spanned by eight eigenvectors associated
with the eigenvalues of smallest magnitude. As for the recycle space, the
invariant subspace is computed for the preconditioned operator. For the primary system,
we use the right invariant subspace. For the dual system, we use the left invariant subspace.
As the recycle space improves, the principal angles tend to zero, and so the cosines
tend to one. Consider the results for the primary system. At the first IRKA step and
the end of the first cycle, we see that the recycle space captures a subspace of
dimension four of the invariant subspace. The recycle space gets more accurate at
the end of the second cycle and captures a subspace of dimension seven.
For the second IRKA step, we have a new shift, and so the matrix changes.
Therefore, at the start of the first cycle, we see a slight deterioration of
the recycle space (almost negligible). This recycle space leads to the dashed
curve in Figure \ref{figure:irka5k}. By the end of the first cycle (at the second IRKA step),
all eight eigenvectors are captured. The results for the dual system recycle space are similar.


For the second set of experiments, we reduce the models to $r=3$ degrees of freedom, using as initial shifts, $1.00 \times 10^{-5}, 7.08 \times 10^{-3}$, and $5.01$. We compute the recycle space at every fifth IRKA step. The results are given in Table \ref{table:time}. We implement recycling for the smallest shift only. The linear systems corresponding to the two (larger) shifts converge fast, so recycling Krylov subspaces is not useful for these.
{\it Total iteration count} refers to the sum of iteration counts for solving linear systems
over all shifts and all IRKA steps. {\it Total time} is the time in seconds required by IRKA to converge to the ideal shifts.
This includes the time for all IRKA computations as well as all linear solves
(BiCG or RBiCG, as the case may be). The table illustrates that computing the reduced model without recycling takes about 50\% more time than with recycling. Obviously, the improvement for just the pairs of linear systems where recycling is actually used
is substantially larger.
\begin{table}[h]
   \footnotesize
    \begin{center}
        \begin{tabular}{cccccccc@{\extracolsep{7pt}}ccc}
            \hline \\ [-1.5ex]
            \multirow{3}{*}{Size} & \multirow{3}{*}{s} & \multirow{3}{*}{k}  & \multirow{3}{*}{Drop} & \multirow{3}{*}{IRKA} & \multicolumn{3}{c}{Total iteration count} & \multicolumn{3}{c}{Total time (s)} \\[5pt]
            \cline{6-8} \cline{9-11} \\ [-1.5ex]
            & & & tol & steps & BiCG & RBiCG & Ratio & BiCG & RBiCG & Ratio \\ [5pt]
            \hline \\ [-1.5ex]
            20209 & 40 & 20 & 0.01 & 31 & 3032 & 1434 & 2.11 & 73.82 & 54.28 & 1.36 \\[5pt]
            79841 & 50 & 20 & 0.005 & 44 & 6324 & 2547 & 2.48 & 742.83 & 505.09 & 1.47 \\[5pt]
            \hline \\ [-1.5ex]
        \end{tabular}
    \end{center}
    \caption{\label{table:time} The total number of iterations and computation time over all IRKA iterations
     with BiCG and with RBiCG for the linear systems with the smallest shift. Total time includes
     the time for all computations to compute the reduced model.}
\end{table}


\section{Conclusion}\label{conclSec}
We focus on efficiently solving sequences of dual linear systems. For several classes of problems,
such as the linear systems arising in interpolatory model reduction, or bilinear forms arising
in Quantum Monte Carlo methods, the BiCG algorithm has advantages over methods like
GMRES that  would solve the primary and the dual system separately.
For sequences of dual linear systems arising in such problems, it is advantageous to use
Krylov subspace recycling for the BiCG algorithm, and for this purpose we propose the RBiCG algorithm.
The derivation of RBiCG also provides the foundation for recycling variants of
other popular bi-Lanczos based methods, like CGS, BiCGSTAB, QMR, and TFQMR \cite{ahuja:MS}.

We have demonstrated the usefulness of RBiCG for interpolatory model reduction using IRKA,
an application that may be an important niche for this solver.
In addition, we have analyzed and demonstrated the effectiveness of RBiCG
for nonsymmetric linear systems arising from convection-diffusion problems.
This suggests that the RBiCG method
may be useful in other areas
where solving dual systems in a Petrov-Galerkin sense brings special advantages.

In future work, we plan to extend the use of RBiCG to model reduction for MIMO dynamical systems
in a tangential interpolation framework where the right-hand sides are not constant as in the SISO case.
In addition, we will investigate the use of RBiCG for evaluating bilinear forms arising in QMC algorithms.
Our current results for this look promising.

\vspace{6pt}
\noindent
{\bf Acknowledgments.} We thank the anonymous reviewers for their careful and helpful
suggestions, which greatly helped us to improve this paper.

\bibliographystyle{abbrv}
\bibliography{rbicg}

\end{document}